\documentclass[11pt]{amsart}

\usepackage{amsthm,amscd,amssymb,amsmath,tikz-cd,mathtools,enumerate, ifthen,bm,paralist,marginnote, appendix}
\usepackage{mathrsfs}
\usepackage{subfiles}

\theoremstyle{plain}
\newtheorem{defn}{Definition}[section]
\newtheorem{lemma}[defn]{Lemma}
\newtheorem{theorem}[defn]{Theorem}
\newtheorem{proposition}[defn]{Proposition}
\newtheorem{corollary}[defn]{Corollary}
\newtheorem{conjecture}[defn]{Conjecture}

\theoremstyle{remark}
\newtheorem{remark}[defn]{Remark}
\newtheorem{example}[defn]{Example}

\DeclareOldFontCommand{\bf}{\normalfont\bfseries}{\mathbf}
\newcommand{\HH}{{\mathrm{H}}}
\newcommand{\II}{{I}}

\newcommand{\rank}{{\mathrm {rank} \,}}

\newcommand{\Hom}{{\mathrm {Hom}}}
\newcommand{\Gr}{\mathrm{Gr}}
\newcommand{\Bl}{\mathrm{Bl}}

\newcommand{\hooklongrightarrow}{\lhook\joinrel\longrightarrow}

\usepackage[margin=1in]{geometry}
\begin{document}
\date{}
\author{{Andrew Harder}}
\title{{Mixed Hodge structures in log symplectic geometry}}

	\address{Andrew Harder,
	\textnormal{Department of Mathematics, Lehigh University, Christmas-Saucon Hall, 14 E. Packer Ave., Bethlehem, PA, USA, 18015.}
	\textnormal{\texttt{anh318@lehigh.edu}}}

\maketitle
\begin{abstract}
We study the cohomology rings of snc log symplectic pairs $(X,Y)$ which have log symplectic forms of pure weight.  We show that under a certain natural condition, the cohomology ring of $X \setminus Y$ exhibits the curious hard Lefschetz property. Analogous results are shown to hold for limit mixed Hodge structures associated to good degenerations of projective irreducible holomorphic symplectic manifolds. We provide several examples of log symplectic pairs of pure weight including a class of cluster-type varieties, and examples coming from the work of Feigin and Odesski. We show that the components of the central fiber of good degenerations of projective irreducible holomorphic symplectic manifolds produce log symplectic pairs.
\end{abstract}

\section{Introduction}

The main goal of this paper is to study the interactions between holomorphic symplectic geometry and mixed Hodge structures. The motivation for this comes from the cohomology ring of a log symplectic pair $(X,Y)$.
\begin{defn}
Let $X$ be a smooth projective variety, and let $Y$ be a reduced divisor in $X$. Then $(X,Y)$ is called {\em log symplectic} if there exists a form $\sigma \in \HH^0(X,\Omega^2_X(\log Y))$ so that the map
\[
\sigma(-,-) : \mathscr{T}_X(- \log Y) \times \mathscr{T}_X (-\log Y) \longrightarrow \mathscr{O}_X
 \]
is nondegenerate at every point in $X$. In other words, the map from $\Omega^1_X(\log Y)$ to $\mathscr{T}_X(-\log Y)$ induced by $\sigma$ is an isomorphism of sheaves. If $(X,\emptyset)$ is log symplectic, we call $X$ a {\em projective holomorphic symplectic manifold}. If $Y$ is a simple normal crossings divisor, we call $(X,Y)$ an {\em snc log symplectic pair}.
\end{defn} 
Holomorphic symplectic varieties have been studied extensively, but when $Y$ is nonempty, log symplectic pairs are not as well understood. 

In dimension 2 log symplectic pairs are abundant: let $S$ be a smooth surface and let $D$ be a reduced simple normal crossings divisor on $S$. Then the pair $(S,D)$ is log symplectic. Higher dimensional examples have been constructed by Goto in \cite{goto}. For instance, moduli spaces of $\mathrm{SU}(2)$ monopoles admit log symplectic structures, as do Hilbert schemes of points on any surface $S$ with a fixed anticanonical divisor. Feigin and Odesski constructed Poisson structures on projective spaces coming from vector bundles on elliptic curves (see also work of Polishchuk \cite{pol}). Gualtieri and Pym \cite{gp} showed that in certain cases, these Poisson structures are in fact log symplectic structures. 

Log symplectic pairs can also be found when studying degenerations of irreducible holomorphic symplectic manifolds. Let $\mathscr{X}$ be a K\"ahler manifold and let $\pi : \mathscr{X} \rightarrow \Delta$ be a proper map. We say that $\pi$ is a semistable degeneration if $X_t = \pi^{-1}(t)$ is smooth if $t\neq 0$, $X_0 = \pi^{-1}(0)$ is simple normal crossings, and $\pi$ vanishes at most to order 1 along each component of $X_0$. We say that a semistable degeneration is a {\em good degeneration of holomorphic symplectic varieties} if the relative dimension of $\pi$ is $2d$ and there is some section $\sigma$ of $\Omega_{\mathscr{X}/\Delta}^2(\log X_0)$ so that $\sigma^d$ is nowhere vanishing on $\mathscr{X}$.

The condition that $\sigma^d$ is nowhere vanishing on $\mathscr{X}$ implies that for each $t \neq 0$, $X_t$ is a projective holomorphic symplectic manifold. We could reasonably expect that the components of the singular fiber are log symplectic. It is not difficult to prove the following result.
\begin{theorem}[Theorem \ref{thm:gooddegen}]\label{thm:hkdeg}
Let $\pi: \mathscr{X}\rightarrow \Delta$ be a good degeneration of holomorphic symplectic manifolds, let $W$ be an irreducible component of $X_0$, and let $\partial W$ be the intersection of $W$ with the singular locus of $X_0$. Then $(W,\partial W)$ is an snc log symplectic pair.
\end{theorem}
In this article, we study log symplectic pairs via the the mixed Hodge structure on the cohomology ring $\HH^*(X\setminus Y;\mathbb{Q})$ for $(X,Y)$ an snc log symplectic pair. Our results allow us to place restrictions on the topology of log symplectic pairs. These results should have implications for the geometry and topology of good degenerations of holomorphic symplectic manifolds. Our strongest results hold for log symplectic pairs which are of pure weight. This is an important concept, so we will define it now.

Recall that, by work of Deligne \cite{del-ii,delIII}, there are mixed Hodge structures on the cohomology groups of any variety $U$. This structure consists of a decreasing Hodge filtration $F^\bullet$ on $\HH^*(U;\mathbb{C})$ and an increasing weight filtration $W_\bullet$ on $\HH^*(U;\mathbb{Q})$ satisfying certain compatibility conditions. Deligne defines a functorial decomposition of $\HH^i(U;\mathbb{C})$ into subspaces coming from these two filtrations. This decomposition is called the Deligne splitting.
\[
I^{p,q}(\HH^\ell(U)) = F^p \cap W_{p+q}^\mathbb{C} \cap\left( \overline{F^q} \cap W_{p+q}^\mathbb{C} + \sum_{j \geq 2} \overline{F^{q-j+1}} \cap W_{p+q-j}^\mathbb{C}\right).
\]
Here, $W_\bullet^\mathbb{C}$ denotes the filtration on $\HH^k(U;\mathbb{C})$ induced by $W_\bullet$ after change of coefficients. Given an snc log symplectic pair $(X,Y)$, the logarithmic form $\sigma$ is closed, therefore it defines a class in $\HH^2(X \setminus Y;\mathbb{C})$. We say that $(X,Y)$ has {\em pure weight $w$} if the class of $\sigma$ lies in $I^{2,w}(\HH^2(X\setminus Y))$\footnote{We will show in Proposition \ref{prop:invweight} that the weight $w$ is an invariant of $(X,Y)$ and does not depend on a particular choice of $w$.}. 

Log symplectic pairs of pure weight appear in a number of ways, some of which are discussed in Section \ref{sect:exes}, but they also appear as the components of good degenerations of irreducible holomorphic symplectic manifolds\footnote{Recall that a projective holomorphic symplectic manifold is called {\em irreducible} holomorphic symplectic (IHS) if it is simply connected and if $\oplus_{j=0}^{2d}\HH^0(X;\Omega_X^j)$ is generated as a ring by $\HH^0(X;\mathscr{O}_X)$ and a symplectic form $\sigma$ in $\HH^2(X;\Omega_X^2)$.} (Theorem \ref{thm:gooddegen}).

Log symplectic pairs of pure weight have a number of nice properties. Recall that for any normal crossings divisor $Y$, there is a simplicial complex which encodes the intersections of the irreducible components of $Y$. This is called the dual intersection complex of $Y$.

\begin{theorem}[Theorem \ref{thm:classification}]
If $(X,Y)$ is an snc log symplectic pair with log symplectic form $\sigma$ of pure weight $1$ or $2$ and $\dim X = 2d$, then the dual intersection complex of $Y$ is of dimension $dw-1$. If $\sigma$ has pure weight $0$ then $Y = \emptyset$.
\end{theorem}

Pure weight also places strong restrictions on the mixed Hodge structure of $\HH^*(X\setminus Y ; \mathbb{Q})$. First, it is not hard to see from the definitions (Proposition \ref{prop:symmetry}) that for any snc log symplectic pair $(X,Y)$ with $\dim X = 2d$, we have that
\begin{equation}\label{eq:symm}
\dim \Gr_F^{m} \HH^{\ell}(X\setminus Y;\mathbb{C}) = \dim \Gr_F^{2d-m} \HH^{2d-2m+\ell}(X\setminus Y;\mathbb{C})
\end{equation}
for all $m,\ell$. This is strengthened by Proposition \ref{thm:hodge-filt} which expresses the statement that for an snc log symplectic pair, the Hodge filtration on $\HH^*(X\setminus Y;\mathbb{Q})$ is determined by $\sigma$ and the cup product structure.

In the case where $(X,Y)$ is snc log symplectic and $\sigma$ has pure weight, we also obtain restrictions on the weight filtration on $\HH^*(X\setminus Y;\mathbb{Q})$. An example of this is the following.
\begin{theorem}[Theorem \ref{thm:chl}, Corollary \ref{cor:scr-geom}]\label{thm:ht}
If $(X,Y)$ is an snc log symplectic pair so that $\dim X = 2d$ and $\sigma$ is a log symplectic form pure weight 2, then
\begin{enumerate}
    \item $I^{p,q}(\HH^\ell(X\setminus Y)) = 0$ if $p \neq q$,
    \item $I^{m,m}(\HH^{\ell}(X\setminus Y)) \cong I^{2d-m,2d-m}(\HH^{2d-2m+\ell}(X\setminus Y))$
for all $\ell,m$.
\end{enumerate}

\end{theorem}
A similar, but less striking, result holds in the case of pure weight 1 (Theorem \ref{thm:ehl}, Corollary \ref{cor:scr-geom}), which says that the weight-graded pieces of $\HH^*(X\setminus Y ;\mathbb{C})$ vanish in certain ranges. These results, combined with the fact that for any smooth variety $U$, $\Gr^W_j \HH^i(U;\mathbb{Q}) \cong 0$ if $j < i$, leads to the following result.
\begin{theorem}[Proposition \ref{prop:pw2}, Proposition \ref{prop:vanpw2}]
Let $(X,Y)$ be a log symplectic pair of pure weight $w$.
\begin{enumerate}
    \item If $w=1$ then $\dim \Gr_F^i\HH^j(X\setminus Y;\mathbb{C}) = 0$ if $j-i < d$.
    \item If $w=2$ then $\dim \Gr_F^i\HH^j(X\setminus Y;\mathbb{C}) = 0$ if $j > 2d$ or $i < j/2$.
\end{enumerate}
\end{theorem}
\begin{figure}[ht]
\begin{centering}
\label{fig:hdpw2}
\begin{tikzpicture}[scale=0.7]
\draw (0,0) -- (-3.5,-3.5) -- (0,-7) -- (3.5,-3.5) -- (0,0);
\draw (0,0) -- (0,-3.5) -- (-3.5,-3.5);

\foreach \i in {0,...,6}
{
\foreach \j in {0,...,6}
{
\ifthenelse{ \not{\j > \i} \AND \numexpr \i + \j < 7}
{\node  at (0 +\j/2 - \i/2,-0.5-\j/2-\i/2) {$*$};}{\node  at (0 +\j/2 - \i/2,-0.5-\j/2-\i/2) {$0$};}
}
}
\draw [dashed] (-1.75,-1.75) -- (1.75,-5.25);
\end{tikzpicture}\quad \quad \quad \begin{tikzpicture}[scale = 0.7]
\draw (0,0) -- (-3.5,-3.5) -- (0,-7) -- (3.5,-3.5) -- (0,0);
\draw (-1.75,-5.25) -- (1.75,-1.75);
\foreach \i in {0,...,6}
{
\foreach \j in {0,...,3}
{
\node  at (0 +\j/2 - \i/2,-0.5-\j/2-\i/2) {$*$};
}
}
\foreach \i in {0,...,6}
{
\foreach \j in {0,...,2}
{
\node  at (2 +\j/2 - \i/2,-2.5-\j/2-\i/2) {0};
}
\draw [dashed] (-1.75,-1.75) -- (1.75,-5.25);
}
\end{tikzpicture}
\caption{Vanishing of the Hodge diamond of a log symplectic pair of pure weights 2 and 1 respectively. The dashed lines indicate the axis of symmetry.}
\end{centering}
\end{figure}
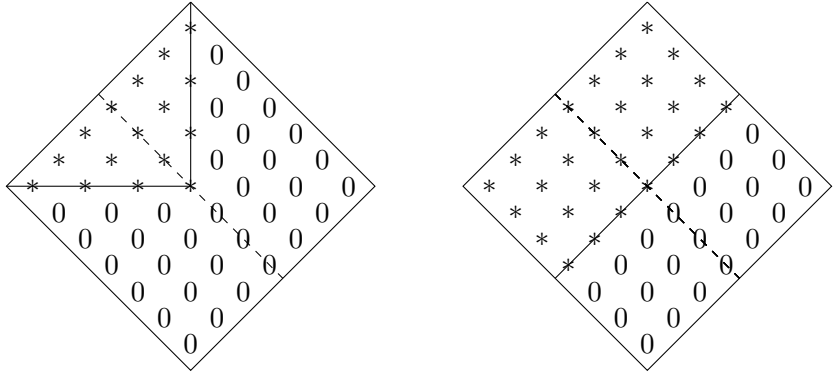

The symmetries in Theorem \ref{thm:ht} also appear in the Hodge theory of character varieties, particularly, the work of Hausel and Rodriguez-Villegas \cite{hrv}. If $C$ is a smooth curve and $G$ is an algebraic group, one can produce the twisted character variety of $C$ associated to $G$. Precisely, we choose a point $p \in C$ and take representations of the fundamental group $\pi_1(C \setminus p)$ in $G$ whose monodromy along a small loop around $p$ is an appropriate multiple of the identity matrix by a root of unity. If we call the collection of all such representations $\mathrm{Rep}(C,G)$, then the twisted character variety $\mathrm{M}_{C,G}$ is the quotient $\mathrm{Rep}(C,G)  /\!\!/  G$ where $G$ acts by conjugation. If $G$ is either $\mathrm{SL}_2(\mathbb{C}), \mathrm{GL}_2(\mathbb{C})$ or $\mathrm{PGL}_2(\mathbb{C})$, then Hausel and Rodriguez-Villegas show that $\HH^*(\mathrm{M}_{C,G};\mathbb{Q})$ has the properties described in Theorem \ref{thm:ht} (see \cite{hrv}). In fact, their results are stronger than this. 
\begin{defn}
A variety $U$ of dimension $d$ has the {\em curious hard Lefschetz property} if there is some $\alpha$ in $\Gr^W_4\HH^2(U;\mathbb{Q})$ so that the maps induced by the cup product pairing,
\[
(\alpha^{(d-m)} \cup (-)) : \Gr^W_{2m}\HH^{\ell}(U;\mathbb{Q}) \longrightarrow \Gr^W_{4d-2m}\HH^{2d-2m+\ell}(U;\mathbb{Q})
\]
are isomorphisms for all $\ell,m$.
\end{defn}
Hausel and Rodriguez-Villegas show that if $G = \mathrm{GL}_2(\mathbb{C}), \mathrm{SL}_2(\mathbb{C})$ or $\mathrm{PSL}_2(\mathbb{C})$ then $\mathrm{M}_{C,G}$ has the curious hard Lefschetz property. Recently, this has also been proven for $\mathrm{GL}_n(\mathbb{C})$ by Mellit \cite{mellit}. One of the main theorems in this paper is the following.
\begin{theorem}[Theorem \ref{thm:chl}, Corollary \ref{cor:scr-geom}]\label{thm:ls2}
Let $(X,Y)$ be a log symplectic pair of pure weight $2$. Then $X \setminus Y$ has the curious hard Lefschetz property.
\end{theorem}
\begin{remark}
In fact the proof of this result applies to more than just log symplectic pairs. For instance, our proof also shows that the limit mixed Hodge structure of a type III degeneration of projective IHS manifolds has the curious hard Lefschetz property (Corollary \ref{cor:chl-deg}). This is related to a result of Soldatenkov \cite{sold}.
\end{remark}

\subsubsection*{Organization}

We will begin in Section \ref{sect:cohomology} by defining a class of objects called symplectic Hodge rings (see Definition \ref{defn:hdgrng}) which formalize the properties of the cohomology rings of log symplectic pairs. We then study the formal properties of symplectic Hodge rings, observing that their Hodge filtration is determined by the algebra structure and the class of $\sigma$ (Proposition \ref{thm:hodge-filt}). We then prove several theorems about the structure of symplectic Hodge rings when $\sigma$ is of pure weight (Theorem \ref{thm:chl}, Theorem \ref{thm:ehl}).

The next three sections provide examples of symplectic Hodge rings. Section \ref{sect:hfilt} focuses on showing that the cohomology rings of snc log symplectic pairs are symplectic Hodge rings. We provide several examples of log symplectic pairs of pure weight 1, from Hilbert schemes of points on surfaces, and from purely elliptic log symplectic pairs. We then generalize a result of Hacking and Keel and show that, given a log symplectic pair $(X,Y)$, one may produce a new log symplectic pair by blowing up the symplectic leaves of the associated Poisson structure on $X$. We also show that this procedure preserves pure weight.

In Section \ref{sect:degenerations}, we show that if $\pi : \mathscr{X}\rightarrow \Delta$ is a good degeneration of compact holomorphic symplectic manifolds, then the limit mixed Hodge ring $\HH^*(X_\infty;\mathbb{Q})$ is a symplectic Hodge ring. We then look at the case where a general fiber $X_t$ of $\pi$ is a projective IHS manifold, in which case $\HH^*(X_\infty;\mathbb{Q})$ has pure weight depending on the order of nilpotency of the monodromy operator on $\HH^2(X_\infty;\mathbb{Q})$. We will make some remarks about the relationship between this result and Nagai's conjecture (Theorem \ref{thm:nag}) in light of the results in Section \ref{sect:cohomology}. Finally, we relate good degenerations and log symplectic pairs: we show that if $\pi : \mathscr{X} \rightarrow \Delta$ is a good degeneration of projective holomorphic symplectic manifolds and $\HH^*(X_\infty;\mathbb{Q})$ has pure weight $w$, then each component of $\pi^{-1}(0)$ is log symplectic of pure weight $w$ as well (Theorem \ref{thm:gooddegen}).

In Appendix \ref{sect:mhc} we prove a technical result necessary for the proof of Theorem \ref{thm:gooddegen} which we were not able to find in the literature, but which is surely known to experts. 
\subsubsection*{Conventions}

In this paper, all varieties are taken over $\mathbb{C}$. Cohomology groups are usually written with their coefficients, but if the coefficient group is absent, the reader should assume that coefficients are complex.

If $V$ is a rational vector space, we will use the notation $V_\mathbb{C}$ to mean $V \otimes \mathbb{C}$. We will often discuss relationships between Hodge and weight filtrations on different vector spaces. To avoid introducing new names for the Hodge and weight filtrations on different vector spaces, we will sometimes use the notation $F^\bullet V_\mathbb{C}$ or $W_\bullet V$ to indicate the Hodge and weight filtrations on $V$ if $V$ is unclear. Similarly, when talking about graded pieces of a filtration, we will use the notation $\Gr^p_F$ and $\Gr^W_i$ when there is no ambiguity as to which vector space $F$ and $W$ belong, and we will use the notation $\Gr^p_FV$ and $\Gr^W_iV$ if $V$ is unclear.

\section{Symplectic Hodge rings}\label{sect:cohomology}

Many of the results in this paper are true for both limit mixed Hodge structures of good degenerations of holomorphic symplectic varieties (Section \ref{sect:degenerations}) and for log symplectic pairs (Section \ref{sect:hfilt}). Instead of proving all of our results twice, we will formalize the properties that the cohomology ring of a log symplectic pair and the cohomology ring of a good degeneration share under the name ``symplectic Hodge rings'' and prove general properties about these objects. There are other contexts in which we expect that symplectic Hodge rings appear but are not discussed in this paper. For instance, finite symplectic quotients of log symplectic pairs and their symplectic resolutions, and certain degenerations of holomorphic symplectic manifolds with certain singularities.

\subsection{Mixed Hodge structures}
From this point on, we assume basic familiarity with the formalism of mixed Hodge structures. A good introduction is \cite[Chapter 3]{ps}. The following facts will often be used without remark in the remainder of this section. Let $(V,F,W)$ be a mixed Hodge structure, then one may define Deligne's splitting of $V$ in the following way
\[
{I}^{p,q}(V) = F^p \cap W_{p+q}^\mathbb{C} \cap\left( \overline{F^q} \cap W_{p+q}^\mathbb{C} + \sum_{j \geq 2} \overline{F^{q-j+1}} \cap W_{p+q-j}^\mathbb{C}\right).
\]
This plays the same role in mixed Hodge theory as the splitting $H^{p,q} = F^p \cap \overline{F^q}$ plays for pure Hodge structures.

\begin{proposition}\label{prop:Ipqprops}
Deligne's splitting has the following properties.
\begin{compactenum}
\item $W_k^\mathbb{C} = \bigoplus_{p+q \leq k} \II^{p,q}(V)$.
\item $F^p = \bigoplus_{r \geq p} \II^{r,s}(V)$.
\item $\II^{p,q}(V) \simeq \overline{\II^{q,p}(V)} \bmod \bigoplus_{r < p, s <q} \II^{r,s}(V)$.
\item Homomorphisms of mixed Hodge structures preserve $\II^{p,q}(V)$.
\end{compactenum}
\end{proposition}
A proof of the first two properties may be found in \cite[Lemma-Definition 3.4]{ps}, and the third and fourth can be found in \cite[Theorem 2.13]{cks}. Furthermore, this splitting is such that the natural map from $\II^{p,q}(V)$ to $\Gr_F^p\Gr^W_{p+q}V$ is an isomorphism. If $f \in \Hom_\mathrm{MHS}(V,V')$, then there are induced maps
\[
\Gr_F^p \Gr^W_i V \longrightarrow \Gr_F^p \Gr^W_i V'
\]
which make the diagram
\begin{equation}\label{eq:comdiagmhs}
\begin{tikzcd}
\II^{p,q}(V) \ar[r] \ar[d] & \II^{p,q}(V') \ar[d]\\
\Gr_F^p \Gr^W_{p+q} V \ar[r] & \Gr_F^p \Gr^W_{p+q} V'
\end{tikzcd}
\end{equation}
commute. 

Recall that if $(V,F,W)$ and $(V', F, W)$  are mixed Hodge structures, then $V \otimes V'$ admits a mixed Hodge structure with the following filtrations.
\begin{align*}
W_k (V \otimes V') &= \sum_{a + b = k} W_aV \otimes W_{b}V' \\
F^p (V_\mathbb{C} \otimes_\mathbb{C} V'_\mathbb{C}) &= \sum_{r+s = p} F^rV_\mathbb{C} \otimes_\mathbb{C} {F}^{s}V'_\mathbb{C}.
\end{align*}
The following two propositions are straightforward consequences of the definitions and linear algebra. We leave their verification to the reader.
\begin{proposition}\label{lem:tensorinj}
If $(V,F,W)$ and $(V',F,W)$ are a pair of mixed Hodge structures, then
\[
\II^{a,b}(V \otimes V') = \bigoplus_{\substack{p+p' =a \\ q + q' = b}} (\II^{p,q}(V) \otimes \II^{p',q'}(V')).
\]
\end{proposition}

\begin{proposition}\label{prop:gradeds}
Let $(V,F,W)$ and $ (V',F,W)$ be a pair of mixed Hodge structures. There is an injection of pure Hodge structures,
\[
\Gr^W_aV \otimes \Gr^W_bV' \hooklongrightarrow \Gr^W_{a+b}(V \otimes V').
\]
\end{proposition}

\subsection{Symplectic Hodge rings}
We will now define the main object in this section.
\begin{defn}\label{defn:hdgrng}
A {\em Hodge ring} is a finite dimensional graded $\mathbb{Q}$-algebra $\bm{H} = \oplus_{i=0}^{4d} H^i$ with $H^0 \cong \mathbb{Q}\cdot  \mathrm{id}$ with the following properties.
\begin{compactenum}
\item Each $H^i$ carries a mixed Hodge structure so that $\Gr_F^{j}H_\mathbb{C}^i \cong 0$ for $j < \max\{0,i-2d\}$ and $j > \min\{i,2d\}$. Consequently, if $i \leq 2d$ then $\dim \Gr^W_jH^i \cong 0$ for $j \geq 2i$ or $j < 0$, and if $i \geq 2d$ then $\dim \Gr^W_jH^i \cong 0$ for $j \leq i -2d$ or $j\geq i + 2d$.
\item The multiplication maps
\[
H^i \otimes H^j \longrightarrow H^{i+j}
\]
induced by the ring operation are homomorphisms of mixed Hodge structures. 
\end{compactenum}
Let $\bm{H}_\mathbb{C} = \bm{H} \otimes \mathbb{C}$ and let $H^i_\mathbb{C}$ denote $H^i \otimes \mathbb{C}$. Assume that $\bm{H}$ also has the following property.
\begin{compactenum}
\setcounter{enumi}{2}
\item There is a nonzero element $\sigma \in F^2H^2_\mathbb{C}$ so that the multiplication map
\[
\sigma^{d-m} : \Gr_F^{m}H_\mathbb{C}^{\ell} \longrightarrow \Gr_F^{2d-m}H_\mathbb{C}^{2d-2m+\ell}
\]
is an isomorphism for all $0 \leq m \leq d$ and all $\ell \geq 0$.
\end{compactenum}
We call $\sigma$ a {\em symplectic element}. If a Hodge ring $\bm{H}$ has a symplectic element, we call $\bm{H}$ a {\em symplectic Hodge ring}. We say that $\sigma$ has {\em pure weight $w$} if $\sigma \in \II^{2,w}(H^2)$.
\end{defn}
\begin{remark}\label{rmk:hgrng}
\begin{enumerate}
    \item For the sake of simplicity, we will use the notation $\II^{p,q;j}$ to denote $\II^{p,q}(H^j)$ where the underlying Hodge ring is understood.
    \item The cohomology ring of any complex algebraic variety equipped with its cup product is a Hodge ring. Condition 1 is a consequence of, e.g., \cite[Theorem 5.39]{ps}, and Condition 2 follows from Proposition \ref{prop:del} in the next section. The limit mixed Hodge structure associated to a semistable degeneration also produces a Hodge ring. This will be explained in Section \ref{sect:degenerations}.
    \item
    Note that $\Gr_F^jH_\mathbb{C}^{2(d+1)} = 0$ for $j \geq 2d$ by Condition 1, hence $F^{2(d+1)}H^{2(d+1)}_\mathbb{C}$ is trivial. Because $\sigma \in F^2H_\mathbb{C}^2$, and the multiplication map is compatible with the Hodge filtration, $\sigma^{d+1} = 0$, 
\end{enumerate}

\end{remark}


\subsection{The Hodge filtration in a symplectic Hodge ring}

We now explain that we may reconstruct the Hodge filtration on a the graded pieces of a symplectic Hodge ring $\bm{H}$ from $\sigma$ and the product structure on $\bm{H}_\mathbb{C}$. In this section, we do not impose the condition that $\sigma$ is of pure weight at all. 

\begin{defn}\label{defn:primfilt}
Let $\bm{H} = \oplus_{i=0}^{4d}H^i$ be a Hodge ring, and let $\sigma$ be a symplectic form. We define $G^\bullet_\sigma$to be  the following filtration using $\sigma$. Let 
\begin{align*}
G_\sigma^1  = \ker &\left\{ \sigma^d : \bm{H}_\mathbb{C} \longrightarrow \bm{H}_\mathbb{C} \right\} \\ 
G_\sigma^{2d} = \mathrm{im} &\left\{ \sigma^d : \bm{H}_\mathbb{C} \longrightarrow \bm{H}_\mathbb{C} \right\}.
\end{align*}
Then for each $1 \leq m \leq d$, we define $G_{\sigma}^{m+1}$ and $G_\sigma^{2d-m}$ inductively to be 
\begin{align*}
G_\sigma^{m+1}  = \ker &\left\{ \sigma^{(d-m)} : G_{\sigma}^m \longrightarrow \bm{H}_\mathbb{C}/G_\sigma^{2d-m+1}\right\} \\ 
G_\sigma^{2d-m} = \mathrm{im}&\left\{ \sigma^{(d-m)} : G_\sigma^m \longrightarrow \bm{H}_\mathbb{C}\right\} + G_{\sigma}^{2d-m+1}.
\end{align*}
\end{defn}
Given a nilpotent operator $\eta$ on any vector space $V$, one can produce a unique increasing filtration $L^\eta_\bullet$ on $V$ so that $\eta( L^\eta_i) \subseteq L^\eta_{i-2}$ and $\eta$ induces isomorphisms,
\[
\eta^i: \Gr^{L^\eta}_{i} \xrightarrow{\cong} \Gr^{L^\eta}_{-i}.
\]
One may use this property to check that $G_\sigma^\bullet$ is precisely $L^{\sigma}_{2d-\bullet}$. Furthermore, we may equip $\bm{H}_\mathbb{C}$ with the total Hodge filtration,
\[
F^m \bm{H}_\mathbb{C} = \oplus_{ i= 0}^{4d} F^m H^i_\mathbb{C}
\]
and by assumption, for each $m$ we get isomorphisms
\[
\sigma^{d-m} : \Gr_F^m \bm{H}_\mathbb{C} \longrightarrow \Gr_F^{2d-m}\bm{H}_\mathbb{C}
\]
for all $m$. We obtain the following result.

\begin{proposition}\label{thm:hodge-filt}
Let $\bm{H}$ be a Hodge ring and let $\sigma$ be an element in $F^2H^2$. The filtration $G_\sigma^\bullet$ is equal to the Hodge filtration on $\bm{H}$ if and only if $\sigma$ is a symplectic element.
\end{proposition}

\begin{remark}
Proposition \ref{thm:hodge-filt} may be summarized by saying that once the degree grading on $\bm{H}$ and the weight filtration are fixed, the Hodge structure on $\bm{H}$ is completely determined by the ring structure and by $\sigma$. For example, if $X$ is a projective IHS manifold then we will see that its cohomology ring $\HH^*(X;\mathbb{Q})$ with the symplectic element given by a holomorphic symplectic form is a symplectic Hodge ring. Then $\sigma$ in fact spans $F^2\HH^2(X;\mathbb{C})$ and determines the Hodge filtration on $\HH^2(X;\mathbb{C})$ by letting
\[
F^1\HH^2(X;\mathbb{C}) = \left\{ a \in \HH^2(X;\mathbb{C}) : q(\sigma,a) = 0\right\}, \quad F^0\HH^2(X;\mathbb{C}) = \HH^2(X;\mathbb{C})
\]
where $q$ is the Beauville--Bogomolov--Fujiki form. Then the Hodge structure on $\HH^2(X;\mathbb{Q})$ is equivalent to a choice of $\sigma$ and the weight filtration is trivial. Therefore the Hodge structure on $\HH^2(X;\mathbb{Q})$ determines the Hodge structure on $\HH^i(X;\mathbb{Q})$ for all $i$. 
\end{remark}


\subsection{A key lemma}

In the following sections, we will determine the extent to which the {\em weight} filtration on $\bm{H}$ is determined by $\sigma$. For this, the pure weight condition must be taken into account. The following lemma  demonstrates its importance.

\begin{lemma}\label{prop:mixedis}
Assume that $\bm{H}$ is a symplectic Hodge ring with $\sigma$ a symplectic element of pure weight $w$. Then product with $\sigma^{d-m}$ induces isomorphisms
\[
\II^{m,s;\ell} \longrightarrow \II^{2d-m,s +w(d-m);\ell + 2(d-m)}.
\]
for all  $m \leq d$ and for all $\ell,s$. Consequently, if $\tau \neq 0 \in \II^{m,s;\ell}$ and $n \leq d-m$ then $\tau \cdot \sigma^{n} \neq 0$.
\end{lemma}
\begin{proof}
By Proposition \ref{prop:Ipqprops}, and Definition \ref{defn:hdgrng}(1), if $\tau_1, \in \II^{p_1,q_1;\ell_1}$ and $\tau_2 \in \II^{p_2,q_2;\ell_2}$ then the product $\tau_1\cdot \tau_2$ is in $\II^{p_1+p_2,q_1 + q_2;i_1 + i_2}$. Therefore, for some $\tau \in \II^{p,q;i}$, the product $\sigma^k \cdot \tau$ is a class in $\II^{p +2k, q+ kw;i+2k}$. 

By Defintion \ref{defn:primfilt}(3), $\sigma^{d-m}$ induces an isomorphism  
\[
\Gr_F^{m}H^{\ell}_\mathbb{C}\longrightarrow \Gr_F^{2d-m}H^{2d-2m+\ell}_\mathbb{C},
\]
hence the induced map
\[
\bigoplus_{s} \II^{m,s;\ell} \longrightarrow \bigoplus_{s} \II^{2d-m,s+w(d-m);2d-2m+\ell}
\]
induced by the product with $\sigma^{d-m}$ is also an isomorphism. Therefore, $\sigma^{d-m}$ induces an isomorphism between the individual subspaces $\II^{m,s;\ell}$ and $\II^{2d-m,s + w(d-m);2d-2m+\ell}$ for all $\ell,m$ and $s$.
\end{proof}

\begin{proposition}\label{prop:invweight}
Let $\bm{H} = \oplus_{i=0}^{4d} H^i$ be a Hodge ring and assume that $\sigma_1,$ and $\sigma_2$ are symplectic elements of pure weight $w_1$ and $w_2$ respectively. Then $w_1 = w_2$.
\end{proposition}
\begin{proof}
By Definition \ref{defn:hdgrng}, there is an isomorphism between $\Gr_F^0H_\mathbb{C}^0 \cong \mathbb{C}$ and $\Gr_F^{2d}H_\mathbb{C}^{2d} \cong \mathbb{C}$. Furthermore, by Lemma \ref{prop:mixedis}, $\II^{0,0;0} \cong \mathbb{C}$ and $\II^{2d,2w_1;2d} \cong \II^{2d,2w_2;2d}$. By Proposition \ref{prop:mixedis}, this can only happen if $w_1 = w_2$.  
\end{proof}
Therefore the phrase ``symplectic Hodge ring of pure weight $w$'' may be used unambiguously to mean a Hodge ring with symplectic element of pure weight $w$.

\subsection{Symplectic Hodge rings of pure weight 1}
We will now describe some restrictions on the weights and degrees for symplectic Hodge rings with symplectic form of pure weight 1. We begin by producing lower bounds on certain Deligne--Hodge numbers of pure weight 1 symplectic Hodge rings.
\begin{proposition}\label{prop:bounds}
Let $\bm{H}$ be a symplectic Hodge ring with $\sigma$ a symplectic element of pure weight 1. Then $\II^{j,k;2i} \geq 1$ for all $j + k = 3i$, $j,k \geq i$ and $i \leq 2d$.
\end{proposition}
\begin{proof}
Let us look at the maps 
\[
\left(\Gr^W_3H^2\right)^{\otimes i} \longrightarrow \Gr_{3i}^WH^{2i}
\]
of pure Hodge structures. Let $[\sigma]$ be the class of $\sigma$ in $H^{2,1}(\Gr_3^WH^2)$ and let $[\overline{\sigma}]$ be its complex conjugate class in $H^{1,2}(\Gr_3^WH^2_\mathbb{C})$. Let $\overline{\sigma}$ be a lift of $[\overline{\sigma}]$ to $\II^{1,2;2}$. We would like to show that $\sigma^{ k} \cdot \overline{\sigma}^{\ell} \in I^{2k + \ell, \ell + 2k;2(k+\ell)}$ maps to a nonzero class for all $k,\ell$ so that $k + \ell \leq d$. We know that $\sigma^{ \ell} \neq 0$ for all $\ell \leq d$, therefore $\overline{\sigma}^{ \ell} \neq 0$ for all $\ell \leq d$. Since $\overline{\sigma}^{\ell} \in \II^{\ell,2\ell;2\ell}$ and $\ell \leq d$ by assumption, we may apply Lemma \ref{prop:mixedis} to see that $\sigma^{k} \cdot \overline{\sigma}^{\ell}$ maps to a nonzero element of $\II^{2k+\ell,k + 2\ell;2(\ell+k)}$ for all $\ell + k \leq d$ by Lemma \ref{prop:mixedis}.
\end{proof}

\begin{corollary}\label{corollary:injell}
Let $\bm{H}$ be a symplectic Hodge ring so that $\Gr^W_3H^2 \cong \HH^1(E;\mathbb{Q})(-1)$ for an elliptic curve $E$. Then for all $i \leq d$, there is an injective map from $\mathrm{Sym}^i(\HH^1(E;\mathbb{Q})(-1))$ to $\Gr^W_{3i}H^{2i}$.
\end{corollary}
\begin{proof}
Let our notation be as in the proof of Proposition \ref{prop:bounds}. Under our assumptions and the fact that cup product is a symmetric map, the induced maps
\[
\mathrm{Sym}^i\left(\II^{2,1;2} \oplus \II^{1,2;2}\right) \longrightarrow \oplus_{j+k = 3i} \II^{j,k;2i}
\]
are injective since $\sigma$ and $\overline{\sigma}$ span $\II^{2,1;2}$ and $\II^{1,2;2}$ respectively. By (\ref{eq:comdiagmhs}) it follows then that the map
\[
\mathrm{Sym}^i\left(\Gr_F^2\Gr^{W^\mathbb{C}}_3H_\mathbb{C}^2 \oplus \Gr_F^1\Gr^{W^\mathbb{C}}_3H_\mathbb{C}^2\right) \longrightarrow \oplus_{j \geq i} \Gr_F^{j}\Gr^{W^\mathbb{C}}_{3i}H_\mathbb{C}^{2i}
\]
is injective and therefore that 
\[
\mathrm{Sym}^i\left(\Gr^W_3H^2\right) \longrightarrow \Gr^W_{3i}H^{2i}
\]
is injective as required.
\end{proof}
\begin{remark}
The conditions of Corollary \ref{corollary:injell} are satisfied, for instance, by the symplectic Hodge rings associated to limit mixed Hodge structures of good degenerations of projective IHS manifolds of type II (see Section \ref{sect:degenerations} for notation). Other examples of this type are obtained from the purely elliptic log symplectic structures studied by Pym \cite{pym}.
\end{remark}
Condition (2) in Definition \ref{defn:hdgrng} implies that we may have $\Gr^W_jH^i \neq 0$ for any $j$ between 0 and $2i$ if $i \leq 2d$. The following Theorem says that if $\bm{H}$ is a symplectic Hodge ring of pure weight 1, then only some of these weights can appear.
\begin{theorem}\label{thm:ehl}
Assume that $\bm{H} = \oplus_{\ell=0}^{4d}H^\ell $ is a symplectic Hodge ring of pure weight 1, then the following statements are true.
\begin{enumerate}
    \item If $d \leq \ell \leq 2d$ and $j \geq d + \ell$ then $\Gr_j^WH^\ell = 0$ .
    \item Assume that $\dim \Gr_F^{\ell} H_\mathbb{C}^\ell = 1$ if $\ell \leq 2d$ is even and $0$ otherwise. Then $\Gr_j^WH^\ell \cong 0$ if $d \leq \ell < 2d$ and $j \geq d+\ell -1$.
\end{enumerate} 
\end{theorem}
\begin{proof}
First, we will prove (1). Choose some $\II^{p,q;\ell}$ with $d \leq \ell \leq 2d$ and with $p+q > d + \ell$. By Proposition \ref{prop:mixedis}, $\II^{p,q;\ell} \cong \II^{2d-p,q + d-p,\ell + 2(d-p)}$. We know that $\II^{2d-p,q+d-p;\ell+2(d-p)} \cong \II^{q + d-p,2d-p;\ell + 2(d-p)}$ by Hodge symmetry. Then the fact that $p + q > d + \ell$ is is equivalent to the fact that $q + d-p > \ell + 2(d-p)$. By Definition \ref{defn:hdgrng}(1), we know that $\Gr_F^jH^i = 0$ if $j \geq i$, and by Proposition \ref{prop:Ipqprops}, we have that $\Gr_F^jH^i \cong \oplus_{s} \II^{j,s;i}$. Therefore, $\II^{q+d-p,2d-p;\ell+2(d-p)} =0$, and it follows that $\II^{p,q;\ell} = 0$. Since, by Proposition \ref{prop:Ipqprops}, $(\Gr^W_{j}H_\mathbb{C}^\ell)\otimes \mathbb{C} \cong \oplus_{p+q = j}\II^{p,q;\ell}$, it follows that $\Gr^W_{j}H_\mathbb{C}^{\ell} \cong 0$ if $d \leq \ell \leq 2d$ and $j > d +\ell$.

Before we prove (2), we will show that, under the conditions listed in (2), $\II^{\ell,s;\ell} = 0$ for $s \neq \ell/2$. We have that $\dim \II^{2m,m;2m} \geq 1$ by Proposition \ref{prop:bounds}, and by Proposition \ref{prop:Ipqprops}, $\dim \Gr_F^\ell H_\mathbb{C}^\ell \cong \oplus_{j} \II^{\ell,j;\ell}$. Therefore, our assumption that $\dim \Gr_F^{\ell}H_\mathbb{C}^\ell = 1$ if $\ell$ is even and $0$ if $\ell$ is odd means that $\II^{\ell,s;\ell} \cong 0$ if either $\ell$ is odd (and $s$ is arbitrary), or if $\ell$ is even and $s \neq \ell/2$. 

Now let us prove (2). Note that, by Proposition \ref{prop:Ipqprops}, this is equivalent to showing that $\dim \II^{p,q;\ell} \cong 0$ if $p + q = j \geq d + \ell -1$. By the same reasoning as in the proof of (1), we have that $\II^{p,q;\ell} \cong \II^{q+d-p,2d-p;\ell+2(d-p)}$. As in the proof of (1), we must have that $\II^{q+d-p,2d-p;\ell+2(d-p)} \cong 0$ if $p + q > d + \ell$, so it will be enough to show that $\II^{q + d-p,2d-p;\ell +2(d-p)} \cong 0$ if $p + q = d+\ell$, or equivalently, if $q + d - p = \ell + 2(d-p)$. Under this condition, the argument in the second paragraph shows that $\II^{\ell + 2(d-p),2d-p;\ell + 2(d-p)} \neq 0$ if and only if $2d-p = (\ell + 2(d-p))/2$, but this is equivalent to $2d = \ell$. Therefore, $\II^{p,q;\ell} \cong 0$ under the conditions of (2).
\end{proof}

\begin{remark}
Nagai has conjectured that for certain degenerations of projective IHS manifolds, weight-graded pieces of the limit mixed Hodge structure vanish in certain degrees. We will discuss the relationship between Theorem \ref{thm:ehl} and Nagai's conjecture in Section \ref{sect:nagcon}.
\end{remark}


\subsection{The curious hard Lefschetz property}
In Section \ref{sect:pw2}, we will show that if $\bm{H}$ is symplectic Hodge ring of pure weight 2, then $\bm{H}$ is {Hodge--Tate} and has the {curious hard Lefschetz property}. 
\begin{defn}
We say that a a Hodge ring $\bm{H}$  is {\em Hodge--Tate} if $\II^{p,q;j} = 0$ if $p \neq q$.
\end{defn}
If $\bm{H}$ is Hodge--Tate then it has the property that $\Gr_i^WH^j = 0$ if $i$ is odd, by an application of Proposition \ref{prop:Ipqprops}. Furthermore, the condition that $\bm{H}$ is Hodge--Tate implies that $\dim \Gr^W_{2\ell}H^i = \dim \Gr_F^{\ell}H^i_\mathbb{C}$ for all $i$ and $\ell$.
\begin{remark}
If a symplectic Hodge ring is Hodge--Tate, then any symplectic element $\sigma$ must lie in $\II^{2,2;2}$, therefore it is of pure weight 2. 
\end{remark}
\begin{defn}
We say that a Hodge ring $\bm{H} = \oplus_{i=0}^{4d}H^i$ has the {\em curious hard Lefschetz property} if $\bm{H}$ is Hodge--Tate and there is some class $\alpha \in \Gr_4^WH^2$ so that the maps
\begin{equation}\label{eq:chl}
(\alpha \cup (-))^{(d-m)} : \Gr^W_{2m}H^{\ell} \longrightarrow \Gr^W_{4d-2m}H^{\ell + 2(d-m)}.
\end{equation}
is an isomorphism for all $i$ and $j$.
\end{defn}
Let $W^\mathbb{C}_\bullet$ denote the induced filtration on $H^*_\mathbb{C}$. We first show that having the curious hard Lefschetz property over $\mathbb{C}$ implies that the curious hard Lefschetz property over $\mathbb{Q}$.
\begin{proposition}\label{prop:chl}
Assume that $\bm{H}$ is Hodge--Tate and that there is some $\beta \in \Gr^{W^{\mathbb{C}}}_4H^2_\mathbb{C}$ that has the property that
\begin{equation}\label{eq:cpxchl}
(\beta \cup (-))^{(d-m)} : \Gr^{W^\mathbb{C}}_{2m}H^{\ell}_\mathbb{C} \longrightarrow \Gr^{W^\mathbb{C}}_{4d-2m}H^{\ell + 2(d-m)}_\mathbb{C}
\end{equation}
is an isomorphism for all $i$ and $j$. Then $\bm{H}$ has the curious hard Lefschetz property.
\end{proposition}
\begin{proof}
By assumption, we have that $\Gr^W_4H^2 \cong \mathbb{Q}(-2)^{\oplus m}$ for some $m$, so that we get a map for each $\ell,m$,
\[
(\mathbb{Q}(-2)^{\oplus m})^{\otimes (d-m)} \otimes\Gr_{2m}^WH^{\ell}\longrightarrow \Gr_{4d-2m}^WH^{\ell + 2(d-m)}
\]
and therefore, for all $\ell$ and $m$, we get a morphism of mixed Hodge structures
\begin{equation}\label{eq:extend}
\Gr^W_4H^2\longrightarrow \mathrm{Hom}_{\mathrm{MHS}}(\Gr^W_{2m} H^{\ell}, \Gr^W_{4d-2m}H^{\ell + 2(d-m)}).
\end{equation}
Our goal is to show that there is some $\alpha$ so that the image of $\alpha$ is an isomorphism for all $i,j$. By assumption, this is true if we work over the complex numbers. For a pair of complex vector spaces, the subset of $\Hom(V,W)$ made up of isomorphisms is Zariski open, hence there is a nonempty Zariski open subset of $\Gr^{W^\mathbb{C}}_4H^2$ made up of classes $\beta$ for which (\ref{eq:cpxchl}) holds for all $i$ and $j$.  

The Zariski closure of $\mathbb{Q}^m$ in $\mathbb{C}^m$ is the entirety of $\mathbb{C}^m$, therefore, this nonempty Zariski open subset intersects $\Gr^W_4H^2 \subset \Gr^{W^\mathbb{C}}_4H_\mathbb{C}^2$ nontrivially in a Zariski open subset. Since nonempty Zariski open subsets $\mathbb{Q}^m$ contain rational points, this proves the proposition.
\end{proof}
\begin{corollary}\label{cor:chl}
A symplectic Hodge ring is Hodge--Tate if and only if it has the curious Hard Lefschetz property.
\end{corollary}
\begin{proof}
By definition, if $\bm{H}$ has the curious hard Lefschetz property, it is Hodge--Tate, so we prove the converse.

If $\bm{H}$ is Hodge--Tate then $\II^{m,m;\ell} \cong \Gr_{2m}^{W^\mathbb{C}}H^\ell_\mathbb{C}$ and $\II^{m,m,\ell} \cong \Gr^{m}_FH^{\ell}$ for all $\ell,m$. Proposition \ref{prop:mixedis} shows that $\sigma^{i}$ induces an isomorphism between $\II^{m,m;\ell}$ and $\II^{2d-m,2d-m;\ell+2(d-m)}$, which then implies that it induces an isomorphism between $\Gr^{W^\mathbb{C}}_{2m}H^{\ell}_\mathbb{C}$ and $\Gr^{W^\mathbb{C}}_{4d-2m}H^{\ell + 2(d-m)}_\mathbb{C}$. Then  Proposition \ref{prop:chl} allows us to deduce that $\bm{H}$ has the curious hard Lefschetz property.
\end{proof}


\subsection{Symplectic Hodge rings of pure weight 2}\label{sect:pw2}
In the previous section, we remarked that if a symplectic Hodge ring is Hodge--Tate (or, equivalently, if it has the curious Hard Lefschetz property), then it has pure weight 2. In this section, we will prove the converse.
\begin{theorem}\label{thm:chl}
Let $\bm{H} = \oplus_{i=0}^{4d}H^i$ be a Hodge ring. Then $\bm{H}$ is a symplectic Hodge ring of pure weight $2$ if and only if has the curious hard Lefschetz property.
\end{theorem}
\begin{proof}
The main idea in this proof is to use Hodge symmetry and Lemma \ref{prop:mixedis} to derive a contradiction. Let us assume that there is some $\II^{p,q;b} \neq 0$ and $q \neq p$. By Hodge symmetry, $\II^{p,q;b} \cong \II^{q,p;b}$, so there is some $\II^{r,v;b} \neq 0$ so that $r < v$. Furthermore, by Lemma \ref{prop:mixedis}, we know that $\II^{r,v;b} \cong \II^{2d-r,v+2(d-r);b + 2(d-r)}$, and if $r > d$, then $2d-r < d$. Therefore, if there is some $\II^{p,q;b} \neq 0$ with $p \neq q$, then there is some $\II^{m,s;\ell}\neq 0$ with $m < d$ and $m < s$.

Now choose $m,s,\ell$ so that $m < d, m< s$ and so that there is no $\ell',m',s'$ so that $\ell' < \ell$ and $m' \neq s'$ with the property that $\II^{m',s';\ell'} \neq 0$. Our goal now is to derive a contradiction. We use the same argument as in the previous paragraph. As before, our assumptions imply that
\[
\II^{m,s;\ell} \cong \II^{2d-m,s+2(d-m); 2d-2m+\ell} \cong \II^{s+2(d-m), 2d-m;2d+2m +\ell} \cong \II^{2m-s;2m-s-d;2(2m-s-d)+\ell} \neq 0. 
\]
The first and third equalities come from Lemma \ref{prop:mixedis}, and the equality in the middle is Hodge symmetry. We know that $s > m$ and that $d> m$. Therefore $2m - s -d < 0$, and therefore, $2(2m-s-d) + \ell < \ell$. Furthermore, note that $2m-s = 2m-s-d$ implies that $d = 0$, which cannot happen since, by definition, a symplectic Hodge ring must have $d \geq 1$. Therefore, we have derived a contradiction.

Therefore, if $\bm{H}$ has a symplectic element of pure weight 2, then it is Hodge--Tate. The claim in the theorem then follows by Corollary \ref{cor:chl}.
\end{proof}
\begin{corollary}\label{cor:shende}
Let $\bm{H}$ be a Hodge ring. Then $\bm{H}$ has the curious hard Lefschetz property if and only if there is some $\sigma \in I^{2,2;2}$ so that $G_\sigma^\bullet \bm{H} = F^\bullet \bm{H}$.
\end{corollary}
\begin{remark}
Corollary \ref{cor:shende} is known for the cohomology rings of character varieties. To our knowledge, the first explicit appearance of this result for character varieties is \cite[Corollary 3]{shende}.
\end{remark}

\section{Log Symplectic pairs}\label{sect:hfilt}

In this section, we give our first examples of a symplectic Hodge ring: the cohomology ring of a log symplectic pair. In Section \ref{sect:sp-shr} we will show that a log symplectic pair produces a symplectic Hodge ring. In Section \ref{sect:cohvan} we look at the consequences of pure weight for log symplectic pairs, first in terms of the geometry of $Y$, then in terms of vanishing results for cohomology groups.

We note that the results in this section depend on the divisor $Y$ being simple normal crossings. The main obstruction to extending these results to cases of more general divisors is our understanding of the relationship between singularities of divisors and mixed Hodge theory. We speculate that the results of this section hold at least when $Y$ has toroidal singularities, and perhaps more broadly.

\subsection{Basic definitions}
The following types of objects have been studied (in greater generality) by several authors; they were first defined by Goto \cite{goto}. More recently, Gualtieri--Pym \cite{gp}, and Pym \cite{pym} have studied log symplectic pairs, as have Lima--Pereira \cite{lp}, and Ran \cite{ran}. We will focus on log symplectic pairs $(X,Y)$ where $Y$ is a simple normal crossings divisor, but we note that this is not required in many of the works that we have mentioned. First, recall the following definition.
\begin{defn}\label{defn:general}
Let $Y$ be a divisor in $X$. We define $\Omega^i_X(\log Y)$ to be the subsheaf of meromorphic $i$-forms on $X$ with poles along $Y$ so that if $h$ is a function whose vanishing defines $Y$, then $h\tau$ and $h d \tau$ are holomorphic forms.
\end{defn}
In the case where $Y$ has normal crossings divisors, we may characterize meromorphic forms with logarithmic poles locally. Take a small holomorphic polydisc $\bm{D}$ centered at $p$ with coordinates $x_1,\dots, x_d$, so that $Y \cap \bm{D} = V(x_1\dots x_k)$ for some $k \leq d$. Then $\Omega^1_X(\log Y)|_{\bm{D}}$ is the $\mathscr{O}_{\bm{D}}$ span of 
\[
d\log x_1,\dots, d\log x_k, d x_{k+1},\dots, d x_d
\]
and  $\Omega_X^m(\log Y)|_{\bm{D}}= \wedge^m \Omega_X^1(\log Y)|_{\bm{D}}$. Now we may give the main definition of this section.

\begin{defn}
Let $(X,Y)$ be a pair consisting of a smooth projective variety $X$ and a divisor $Y$. We say that $\sigma \in \HH^0(X,\Omega^2_X(\log Y))$ is a log symplectic form if one of the following equivalent conditions holds.
\begin{compactenum}
\item The map $\sigma(-,-) : \mathscr{T}_X(-\log Y) \otimes \mathscr{T}_X(-\log Y) \longrightarrow \mathscr{O}_X$ is nondegenerate at every point.
\item The wedge power $\sigma^{\dim X/2} \in \Gamma( \Omega^{\dim X}_X(\log Y)) \cong \Gamma(\omega_X(Y))$ is a nonvanishing global section.
\end{compactenum}
If $(X,Y)$ admits a log symplectic form, we say that it is a {\em log symplectic pair}. If $Y$ is a simple normal crossings divisor, we call $(X,Y)$ an snc log symplectic pair. 
\end{defn}

\begin{proposition}\label{prop:symmetry}
If $(X,Y)$ forms an snc log symplectic pair with $\dim X = 2d$, and $\sigma$ is its log symplectic form, then the map
\[
\sigma^j : \Omega_X^{d-j}(\log Y) \longrightarrow \Omega^{d+j}_X(\log Y), \quad \tau \longmapsto \sigma^j \wedge \tau
\]
is a sheaf isomorphism. Therefore, $\HH^i(X, \Omega_X^{d-j}(\log Y)) \cong \HH^i(X,\Omega_X^{d+j}(\log Y))$ for all $i$ and $j$.
\end{proposition}
\begin{proof}
The map induced by $\sigma^j$ is certainly a morphism of sheaves, since $\sigma$ is a global holomorphic form. Since $\sigma$ is nondegenerate, it follows that follows that it is an injection on stalks. Since $\rank \Omega_X^j(\log Y) = \rank \Omega_X^{d-j}(\log Y)$ it is an isomorphism on stalks and hence is a sheaf isomorphism.
\end{proof}
\begin{remark}
This implies that the Hodge numbers of a log symplectic pair, which we define to be $h^{p,q}(X\setminus Y) : = \dim \HH^q(X,\Omega_X^p(\log Y))$, have the property that $h^{p,q}(X\setminus Y) = h^{2d-p,q}(X\setminus Y)$.
\end{remark}

\subsection{The mixed Hodge structure on the cohomology of a smooth quasiprojective variety}
Our goal is to show that if $(X,Y)$ is an snc log symplectic pair with log symplectic form $\sigma$ and $\dim X = 2d$ then $\oplus_{\ell=0}^{4d} \HH^\ell(X\setminus Y;\mathbb{Q})$ is a symplectic Hodge ring with symplectic element $[\sigma]$. In order to do this, we will review a little bit of the construction which produces mixed Hodge structures on the cohomology of algebraic varieties. We will give full details in Appendix \ref{sect:mhc}.

Let $X$ be a smooth projective variety, let $Y$ be a snc divisor in $X$, and let $U = X \setminus Y$. Let $j : U \hooklongrightarrow X $ be the open embedding. Then there is a quasiisomorphism of complexes,
$\Omega_X^\bullet(\log Y) \rightarrow Rj_*\Omega_U^\bullet$, hence an isomorphism $\mathbb{H}^\ell(X;\Omega_X^\bullet(\log Y)) \cong \HH^\ell(U;\mathbb{C}) \cong \HH^\ell(U;\mathbb{Q})\otimes \mathbb{C}$ for all $\ell$. Deligne has shown that $\HH^\ell(U;\mathbb{Q})$ admits a mixed Hodge structure, and we may identify $F^p\HH^\ell(U;\mathbb{C})$ and $W^\mathbb{C}_\bullet\HH^\ell(U;\mathbb{C})$ as follows.
\begin{defn}\label{defn:hodgeweight}
Let
\[
W_{i}^\mathbb{C}\Omega_X^j(\log Y) = \begin{cases} 0 & \mathrm{if } \quad m < 0 \\
\Omega^j_X(\log Y) & \mathrm{ if } \quad i \geq j \\
\Omega^{j-i}_X \wedge \Omega^i_X(\log Y) & \mathrm{ if } \quad 0 \leq i \leq j\end{cases}
\]
and
\[
F^p\Omega_X^j(\log Y) = (0 \longrightarrow \dots \longrightarrow 0 \longrightarrow \Omega^p_X(\log Y) \longrightarrow \Omega_X^{p+1}(\log Y) \longrightarrow  \dots).
\]
Then the Hodge and (complexified) weight filtration of Deligne's canonical mixed Hodge structure on $\HH^\ell(U;\mathbb{C})$ is given by
\begin{align*}
F^p\mathbb{H}^\ell(X;\Omega_X^\bullet(\log Y)) &= \mathrm{im}(\mathbb{H}^\ell(X;F^p\Omega_X^\bullet(\log Y)\longrightarrow \mathbb{H}^\ell(X;\Omega_X^\bullet(\log Y)) \\
W^\mathbb{C}_{\ell+i}\mathbb{H}^\ell(X;\Omega_X^\bullet(\log Y)) &= \mathrm{im}(\mathbb{H}^\ell(X;W_{i}^\mathbb{C}\Omega_X^\bullet(\log Y)\longrightarrow \mathbb{H}^\ell(X;\Omega_X^\bullet(\log Y))
\end{align*}
\end{defn}
The cohomology ring of $U$ admits a cup product as well. On the level of de Rham cohomology the cup product is given by taking the wedge product of differential forms; on the level of complexes, this map is produced by taking the map in hypercohomology induced by the morphism
\[
\Omega_X^\bullet(\log Y) \otimes \Omega_X^\bullet(\log Y) \longrightarrow \Omega_X^\bullet(\log Y)
\]
induced by wedge product on logarithmic forms. The following statement is a result of Deligne.
\begin{proposition}[{\cite[Corollaire 8.2.11]{delIII}}]\label{prop:del}
The cup product on rational cohomology induces morphisms of mixed Hodge structures,
\[
\cup : \HH^\ell(U;\mathbb{Q}) \otimes \HH^m(U;\mathbb{Q}) \longrightarrow \HH^{\ell + m}(U;\mathbb{Q}).
\]
\end{proposition}
In other words, the cup product in cohomology on a smooth quasiprojective variety preserves mixed Hodge structures. Combining this with Remark \ref{rmk:hgrng}, we obtain the following consequence.
\begin{corollary}
If $U$ is a smooth quasiprojective variety of dimension $d$, then $\bm{H}_U = \oplus_{\ell=0}^{2d} \HH^\ell(U;\mathbb{Q})$ is a Hodge ring.
\end{corollary}

Global sections of $\Omega_X^i(\log Y)$ can be used to produce elements of $\mathbb{H}^i(X;\Omega_X^\bullet(\log Y))$. In general, if $\mathscr{F}^\bullet$ is a complex of sheaves on $X$, and $\mathcal{C}^\bullet_{Gd}\mathscr{F}^\bullet$ is its Godement resolution. Then $\mathbb{H}^i(X,\mathscr{F}^\bullet) = \HH^i(\Gamma(X,\bm{s}(\mathcal{C}_{Gd}\mathscr{F}^\bullet)))$. Therefore, if $\gamma$ is a closed global section of $\mathscr{F}^p$, its image in $\Gamma(\bm{s}(\mathcal{C}^\bullet_{Gd}\mathscr{F}^\bullet))$ is also closed (by construction of $\mathcal{C}_{Gd}^\bullet \mathscr{F}^\bullet$, see e.g. \cite[\S B.2.1]{ps}), hence represents a cohomology class in $\mathbb{H}^p(X,\mathscr{F}^\bullet)$. By the degeneration of the Hodge-to-de Rham spectral sequence for $\Omega_X^\bullet(\log Y)$, any global section $\omega$ of $\Omega_X^i(\log Y)$ is $d$-closed. We denote its class in $\mathbb{H}^i(X;\Omega_X^\bullet(\log Y))$ by $[\omega]$. The classes $[\omega]$ associated to global logarithmic forms span $F^i\mathbb{H}^i(X;\Omega_X^\bullet(\log Y))$. Every such $\omega$ determines a morphism of complexes, $L_\omega:  \Omega_X^\bullet(\log Y) \rightarrow \Omega^{\bullet}_X(\log Y)[-i]$ given by wedge product with $\omega$, and the induced map sends a hypercohomology class $\alpha \in \mathbb{H}^\ell(X,\Omega_X^\bullet(\log Y))$ to $[\omega] \cup \alpha$.
\subsection{Log symplectic pairs with log symplectic forms of pure weight}\label{sect:sp-shr}

In this section, we investigate the consequences of the following definition.
\begin{defn}
Let $(X,Y)$ be a snc log symplectic pair and let $\sigma$ be a log symplectic form on $(X,Y)$. We say that $\sigma$ is of {\em pure weight $w$} if its class in $\HH^2(X\setminus Y;\mathbb{C})$ is in $I^{2,w;2}$.
\end{defn}
Our first observation is that if $\sigma$ is of pure weight, then the geometry of $Y$ is heavily constrained. Recall that to any snc divisor $Y$ of dimension $d$ which is a union of irreducible divisors $Y_1,\dots, Y_n$, there is a simplicial complex called its dual intersection complex. This complex is obtained by taking a simplex $\Delta_I$ of dimension $|I|$ for each subset $I \subset \{1,\dots, n\}$ for which $Y_I = \cap_{i\in I}Y_i$ is nonempty, and attaching $\Delta_I$ to $\Delta_J$ if $Y_J$ is contained in $Y_I$. The dimension of this complex is 1 less than the largest number of components of $Y$ that contain a single point. Since $Y$ is normal crossings and $X$ is assumed to have dimension $2d$, the dual intersection complex of $Y$ has dimension at most $2d-1$.

\begin{proposition}\label{thm:classification}
If $(X,Y)$ is a log symplectic pair of pure weight $w=1,2$ and $\dim X = 2d$, then the dual intersection complex of $Y$ is of dimension $dw-1$. If $(X,Y)$ is a log symplectic pair with $\sigma$ a log symplectic form of pure weight 0, then $Y = \emptyset$.
\end{proposition}
\begin{proof}
By our assumption that $\sigma$ be of pure weight, the class $[\sigma]$ is in $I^{2,w}(\HH^2(U))$ for $i = 0,1$ or $2$. By Proposition \ref{prop:del} cup product respects mixed Hodge structure. Since $\sigma^{\otimes d}$ is the cup product of $\sigma$ with itself $d$ times it is contained in $I^{2d,wd}(\HH^{2d}(U))$ by Proposition \ref{prop:Ipqprops}. Therefore, $\sigma^{\otimes d}$ is a nondegenerate holomorphic form in $\Omega_X^{2d-dw} \wedge \Omega_X^{dw}(\log Y)$ which is not contained in $\Omega_X^{2d-dw + 1} \wedge \Omega_X^{dw-1}(\log Y)$ so we must have $\Omega_X^{2d-dw} \wedge \Omega_X^{dw}(\log Y) \neq \Omega_X^{2d-dw+1} \wedge \Omega_X^{dw-1}(\log Y)$.

We now remark that if $D$ is a polydisc around a point $p$ in complex variables $x_1,\dots, x_{2d}$ in $X$ so that $Y$ is the vanishing locus of $x_1,\dots, x_\ell$ then $\Omega_X^m(\log Y)$ is locally a free module over $\mathscr{O}_D$ whose basis is given by wedge products of 
\[
d \log x_1,\dots, d\log x_\ell, d x_{\ell + 1}, \dots, d x_{2d}.
\]
Therefore, the dimension of the dual intersection complex of $Y$ is the largest $j$ so that $\Omega_X^{2d-j} \wedge \Omega_X^{j}(\log Y) \neq \Omega_X^{2d}(\log Y) $. If $w = 2$, the argument in the previous paragraph shows that $\Omega^{2d}_X(\log Y) \neq \Omega_X \wedge \Omega^{2d-1}_X(\log Y)$, so the dual intersection complex of $Y$ must have dimension at least $2d-1$. Since the dimension of the dual intersection complex of $Y$ is at most $2d-1$, it must be exactly $2d-1$ if $w=2$.  

In the remaining cases, the fact that $\Omega^{2d-dw}_X\wedge\Omega_X^{dw}(\log Y) \neq \Omega_X^{2d-dw+1}\wedge \Omega_X^{dw-1}(\log Y)$ shows the dual intersection complex of $Y$ is of dimension {\em at least} $dw-1$. We now use the fact that $\sigma$ is a non-degenerate form to show that we have equality. Let $\sigma$ be a global section of $\Omega_X^{2d-dw} \wedge \Omega_X^{dw}(\log Y)$. Assume that we may find a local chart with coordinates $(x_1,\dots, x_{2d})$ so that locally $Y = V(x_1\dots x_{dw + j})$ for some $j \geq 0$. Then
\[
\sigma = h(x_1,\dots, x_{2d}) d\log x_1 \wedge \dots \wedge d\log x_{di +j} \wedge d x_{di+j+1} \wedge \dots \wedge dx_{2d}.
\]
Since $\sigma$ is in $\Omega^{2d-dw}_X \wedge \Omega^{dw}_X(\log Y)$, it follows that $h(x_1,\dots,x_{2d})$ must vanish along some subset of the divisors cut out by $x_1,\dots, x_{dw+j}$, hence $\sigma$ is degenerate unless $j = 0$. Therefore, the dimension of the dual intersection complex of $Y$ is at most $dw$. This completes the proof of the proposition.
\end{proof}
\begin{corollary}
If $(X,Y)$ is a log symplectic pair of pure weight 0 then $Y = \emptyset$ and $X$ his holomorphic symplectic.
\end{corollary}
\begin{remark}
One may show via examples that the consequences of Theorem \ref{thm:classification} can fail when $(X,Y)$ is not of pure weight. Let $S_2$ be a K3 surface and let $S_2$ be a rational surface with $E$ a smooth anticanonical elliptic curve in $S_2$. Then $(S_1 \times S_2, S_1 \times E)$ is a log symplectic pair not of pure weight whose boundary divisor has dual intersection complex of dimension 0. 
\end{remark}

We will now  strengthen Proposition \ref{prop:symmetry} to a statement about how cup product with $\sigma$ behaves with respect to the mixed Hodge structure on the Hodge ring $\HH^*(X\setminus Y;\mathbb{Q}) := \oplus_{i=0}^{4d}\HH^i(X\setminus Y;\mathbb{Q})$.

\begin{theorem}\label{eq:filtis}
Let $(X,Y)$ be a snc log symplectic pair and let $\sigma$ be a log symplectic form on $(X,Y)$. The map $L_\sigma^{(i+j)}$ induces an isomorphism
\[
\Gr_F^{d-i-j}\mathbb{H}^{d-i}(X,\Omega_X^\bullet (\log Y)) \longrightarrow \Gr_F^{d+i+j}\mathbb{H}^{d+i+2j}(X,\Omega_X^\bullet (\log Y))
\]
for all $i$ and $j$. Therefore, if $(X,Y)$ is an snc log symplectic pair, then $\HH^*(X\setminus Y;\mathbb{Q})$ is a symplectic Hodge ring. If $\sigma$ is of pure weight $w$ then $\HH^*(X\setminus Y;\mathbb{Q})$ is a symplectic Hodge ring of pure weight $w$.
\end{theorem}
\begin{proof}
Wedge product with $\sigma$ induces a map
\[
F^p\Omega_X^{j}(\log Y)  \longrightarrow  F^{p+2}\Omega_X^{j+2}(\log Y)
\]
by definition of the Hodge filtration. The Hodge to de Rham spectral sequence for the Hodge filtration on $\Omega^\bullet_X(\log Y)$ is the spectral sequence induced by the inclusions 
\begin{equation}\label{eq:degenerate}
F^{p+1} \Omega^\bullet_X(\log Y)  \hooklongrightarrow F^{p} \Omega^\bullet_X(\log Y).
\end{equation}
This spectral sequence degenerates at the $E_1$ term which implies that the maps
\[
\mathbb{H}^k(X, F^p \Omega_X^\bullet(\log Y)) \longrightarrow \mathbb{H}^k(X, \Omega_X^\bullet(\log Y))
\]
are injections. There is a long exact sequence obtained from (\ref{eq:degenerate}) coming from the short exact sequence of complexes,
\[
0\longrightarrow F^{p+1} \Omega^\bullet_X(\log Y)  \longrightarrow F^{p} \Omega^\bullet_X(\log Y) \longrightarrow \Omega^p_X(\log Y)[p] \longrightarrow 0.
\]
Since the corresponding spectral seqence degenerates at the $E_1$ term, it breaks up into a number of short exact sequences of the form
\[
0 \longrightarrow \mathbb{H}^k(X,F^{p+1}\Omega^\bullet_X(\log Y))  \longrightarrow \mathbb{H}^k(X,F^p \Omega^\bullet_X(\log Y)) \longrightarrow \mathbb{H}^{k-p}(X,\Omega_{X}^{p}(\log Y)) \longrightarrow 0
\]
where $\mathbb{H}^k(X,F^p \Omega^\bullet_X(\log Y))$ is isomorphic to $F^p\HH^k(U;\mathbb{C})$. The map $L_\sigma^{(i+j)}$ induces a morphism of filtered complexes $\Omega_X^\bullet(\log Y) \rightarrow \Omega_X^{\bullet + 2k}(\log Y)$ sending $F^p\Omega^\bullet_X(\log Y)$ to $F^{p+2k}\Omega^{\bullet+2k}_X(\log Y)$. The same map sends $\Omega_X^p(\log Y)[p]$ to $\Omega_X^{p+2k}(\log Y)[p]$. If $p = d-k$ then this is a quasiisomorphism of complexes by Proposition \ref{prop:symmetry}. Hence we obtain a commutative diagram of short exact sequences induced by $\sigma^{(i+j)}$,
\[
\begin{tikzcd}
F^{d-i-j+1} \HH^{d-i}(X\setminus Y;\mathbb{C}) \ar[r] \ar[d] & F^{d-i-j} \HH^{d-i}(X\setminus Y;\mathbb{C}) \ar[r] \ar[d] & \HH^{j}(X,\Omega_{X}^{d-i-j}(\log Y)) \ar[d] \\ 
F^{d+i+j+1} \HH^{d+i+2j}(X\setminus Y;\mathbb{C}) \ar[r] & F^{d+i+j} \HH^{d+i+2j}(X\setminus Y;\mathbb{C}) \ar[r] &  \HH^{j}(X,\Omega_{X}^{d+i+j}(\log Y)) 
\end{tikzcd}
\]
which is an isomorphism in the last term. Therefore, we have that $\sigma^{(i+j)}$ induces an isomorphism between $\Gr_F^{d-i-j}\HH^{d-i}(U;\mathbb{C})$ and $\Gr_F^{d+i+j}\HH^{d+i+2j}(U;\mathbb{C})$ for all $i,j,k$.
\end{proof}
The following is then an immediate corollary of Theorem \ref{thm:chl} and Theorem \ref{eq:filtis}.
\begin{corollary}\label{cor:scr-geom}
If $(X,Y)$ is a log symplectic pair with log symplectic form $\sigma$ of pure weight 2 then $\HH^*(X\setminus Y;\mathbb{Q})$ has the curious hard Lefschetz property.
\end{corollary}

\subsection{Cohomological vanishing for log symplectic pairs of pure weight }\label{sect:cohvan}

Now we may discuss the consequences of Theorems \ref{thm:ehl} and \ref{thm:chl} for log symplectic pairs. The new ingredient is the well known fact that if $U$ is a smooth quasiprojective variety, $W_iH^j(U;\mathbb{Q}) = 0$ if $i < j$ (see e.g. \cite[Proposition 4.20]{ps}). This enforces very strict bounds on the possible Deligne--Hodge numbers.

\begin{proposition}\label{prop:vanpw2}
If $(X,Y)$ is a log symplectic pair of pure weight 2 so that $\dim X = 2d$, then $\Gr_F^m \HH^{\ell}(X\setminus Y;\mathbb{C}) = 0$ if $m < \ell/2$ and $\HH^\ell(X\setminus Y ;\mathbb{C}) \cong 0$ if $\ell > 2d$.
\end{proposition}
\begin{proof}
According to Theorem \ref{thm:chl}, the cohomology ring of $U$ is Hodge--Tate. If $U$ is a smooth variety, then $W_i \HH^\ell(X\setminus Y;\mathbb{Q}) = 0$ if $i  < \ell$. On the other hand, being Hodge--Tate implies that $\dim \Gr^W_{2m}H^\ell(X\setminus Y;\mathbb{Q}) = \dim \Gr_F^m \HH^\ell(X\setminus Y;\mathbb{C})$ for all $j$, therefore $\Gr_F^m \HH^\ell(X\setminus Y;\mathbb{C}) = 0$ if $m < \ell/2$. This proves the first statement. 

Now we prove the second statement. Assume that $2m \geq \ell > 2d$. Since $\HH^*(X\setminus Y;\mathbb{Q})$ is Hodge--Tate, $\Gr^{W^\mathbb{C}}_{2m}\HH^\ell(X\setminus Y;\mathbb{C})\cong I^{m,m;\ell} \cong I^{2d-m,2d-m;\ell+2(d-m)}$ (the final isomorphism comes from Proposition \ref{prop:mixedis}). Since $\ell > 2d$, it follows that $2d-m < \ell/2 + (d-m)$, so by the previous paragraph, $\Gr^W_{2m}\HH^\ell(X\setminus Y;\mathbb{Q})$ vanishes for all $2m \geq \ell$. Therefore $\HH^\ell(X\setminus Y;\mathbb{Q})$ vanishes.
\end{proof}

\begin{remark}
Note that an affine variety $U$ of dimension $d$ has the property that $ \HH^i(U;\mathbb{Q}) = 0$ for $i >d$ by the hard Lefschetz theorem. It is easy to see that not all log symplectic pairs of pure weight 2 have the property that $X \setminus Y$ is affine, hence Proposition \ref{prop:vanpw2} contains new information. For instance, one may choose a rational elliptic surface $\mathscr{E}$ with fibration map $\pi : \mathscr{E} \rightarrow \mathbb{P}^1$ so that $\pi$ has a fiber $F$ which is a cycle of rational curves. Then $(\mathscr{E},F)$ is a log symplectic pair but $\mathscr{E} \setminus F$ is not affine, in fact its affinization is $\mathbb{A}^1$.
\end{remark}
\begin{proposition}\label{prop:pw2}
Let $(X,Y)$ be a log symplectic pair of pure weight $1$ and assume that $\dim X = 2d$. Then $\Gr_F^p\HH^\ell(X\setminus Y;\mathbb{C}) = 0$ if $p < \ell-d$.
\end{proposition}
\begin{proof}
Since $\dim \Gr_F^p\HH^\ell(X\setminus Y;\mathbb{C}) = \sum_{s} \dim I^{p,s;\ell}$ it is enough to show that $I^{p,s;\ell} = 0$ for $p < \ell - d$. Before proceeding, we also note that since $X$ is smooth, $\Gr^W_m\HH^\ell(X\setminus Y;\mathbb{Q}) = 0$ if $m < \ell$, hence if $I^{p,s;\ell} \neq 0$ then $p + s \geq \ell$. 

Assume that $I^{p,s;\ell} \neq 0$. Then we have that 
\[
I^{p,s;\ell} \cong I^{s,p;\ell} \cong I^{2d-s,p + d-s,\ell + 2(d-s)}
\]
by Proposition \ref{prop:mixedis} and Hodge symmetry. It follows that if $I^{p,s;\ell} \neq 0$, we must have $3d-2s+p \geq \ell +2(d-s)$, which is equivalent to $p \geq \ell - d$.
\end{proof}

\section{Examples of log symplectic pairs}\label{sect:exes}

In this section, our goal is to construct examples of log symplectic pairs. The most basic examples of log symplectic pairs occur in dimension 2. Let $S$ be a surface, and let $D$ be a snc anticanonical divisor, then $\Omega^2_X(\log D) \cong \omega_S(D) \cong  \mathscr{O}_S$, hence any holomorphic global section of $\omega_S(D)$ is nondegenerate. Furthermore, $\dim \HH^{2,0}(S\setminus D;\mathbb{C}) = 1$ and it follows that $(S,D)$ is of pure weight. Therefore, by Theorem \ref{thm:classification}, being of pure weight 1 is equivalent to $Y$ being smooth, and being of pure weight 2 is equivalent to $Y$ having nodes, and pure weight 0 is equivalent to $D$ being empty, that is, $S$ is a Calabi--Yau surface. 

The following sections explain constructions of log symplectic pairs of higher dimensions.
\subsection{Hilbert schemes of points on log symplectic surfaces}
In this section, we will explain that log symplectic varieties of pure weight 1 exist in all dimensions due to a theorem of Ran. It is classically known that the Hilbert scheme of points on a surface is smooth, and furthermore that if $S$ is a surface admitting a log symplectic structure then there is a log symplectic structure on its Hilbert scheme of points in the sense of Definition \ref{defn:general} (see \cite{goto} for details). 

Precisely, if $(S,D)$ is a log symplectic pair where $\dim S = 2$, then the Hilbert scheme $S^{[n]}$ of $n$ points on $S$ admits a morphism, $S^{[n]} \rightarrow \mathrm{Sym}^n(S)$, called the Hilbert--Chow morphism, to the $n$th symmetric power of $S$. The divisor $D$ determines a divisor $D\rq{}$ in $\mathrm{Sym}^n(S)$ as the image of the divisor 
\[
\widetilde{D} = D \times S \times \dots \times S \subseteq S^n
\]
under the quotient map $S^n \rightarrow \mathrm{Sym}^n(S)$. There is then a log symplectic structure on $S^{[n]}$ whose degeneracy divisor is the proper transform of $D'$ under the Hilbert--Chow morphism $S^{[n]}\rightarrow \mathrm{Sym}^n(S)$. This divisor, which we will call $D''$, is not simple normal crossings if $n \neq 1$, so while the pair $(S^{[n]},D'')$ is log symplectic, it is not {\em snc} log symplectic. The Hilbert scheme parametrizes subschemes of $S$ of length $n$, hence we may stratify $S^{[n]}$ by the intersection index of a length $n$ subscheme with $D$. This is called the incidence stratification of $S^{[n]}$ with respect to $D$.

Ran showed in \cite{ran2} that if $D$ is smooth then iterated blow up along the subschemes making up the incidence stratification makes the preimage of the divisor $D''$ simple normal crossings.  Let $S^{[n]}_r$ denote this blow up and let $D_r$ denote the preimage of $D''$ in $S^{[n]}_r$. Ran showed that $D_r$ is an snc divisor and that the pullback of the log symplectic form on $S^{[n]}$ along this resolution has log poles with respect to $D_r$. Therefore $(S^{[n]}_r,D_r)$ is an snc log symplectic pair.
\begin{proposition}
The snc log symplectic pair $(S^{[n]}_r,D_r)$ is of pure weight 1.
\end{proposition}
\begin{proof}
Ran \cite{ran2} shows that the log symplectic form constructe on $S^{[n]}_r$ is P-normal which means that at each point in $D_r$ there is some local coordinate chart centered at $p$ with coordinates $(z_1,\dots, z_n, y_1,\dots, y_n)$ so that $Y = V(x_1\dots x_k)$ for some $k \leq n$ and 
\begin{equation}\label{eq:Pnormal}
\sigma = d\log z_1 \wedge d y_1 + \dots + d\log z_k \wedge dy_k + d z_{k+1} \wedge d y_{k+1} + \dots + d z_n \wedge d y_n.
\end{equation}
Therefore, $[\sigma] \in W_3\HH^2(S^{[n]}_r \setminus D_r;\mathbb{C})$ by the definition of the weight filtration on the complex $\Omega^\bullet_{S^{[n]}}(\log D_r)$.

Now we would like to show that $I^{2,0;2} = 0$, or equivalently, that $\Gr^W_2\HH^2(S^{[n]}_r \setminus D_r;\mathbb{Q}) \cong \mathbb{Q}(-2)^m$ for some $m$. If we can show this, then since $\sigma \in I^{2,1;2} \oplus I^{2,0;2}$ (by the argument above), it follows that $\sigma\in I^{2,1;2}$. Recall that for any noncompact variety $U$ with a snc compactification $X$, $\Gr_i^W\HH^i(U;\mathbb{Q})$ is the image of the pullback map $\HH^i(X;\mathbb{Q}) \rightarrow \HH^i(U;\mathbb{Q})$ (see e.g. \cite[Remark 8.37]{voisin}). We have a normal crossings compactification $S^{[n]}_r$ of $S^{[n]}_r \setminus D_r$, and furthermore we know that $S^{[n]}$ is unirational. Since $h^{2}(\mathscr{O}_X) = 0$ for any unirational variety, the result follows.
\end{proof}

\subsection{Purely elliptic log symplectic pairs}\label{sect:purel}

We would now like to define what Pym calls purely elliptic log symplectic structures. Quoting the definition in \cite{pym} would take us somewhat far afield, so instead, we choose to provide an an equivalent definition in the case where $\dim X = 4$ coming from the local description in \cite[Theorem 4.5]{pym}.

\begin{defn}[{\cite[Theorem 4.5]{pym}}]
Let $X$ be a smooth fourfold, and let $Y$ be a divisor in $X$. A log symplectic structure $\sigma$ on $(X,Y)$ is {\em purely elliptic} if there are local holomorphic coordinates around each singular point of $Y$ in which $Y$ is given by one of the equations in Table \ref{table:nonlin} and $\sigma$ is of the form
\[
\sigma = d \log f \wedge d w  - a\dfrac{x dy \wedge d z}{\lambda f} + b\dfrac{y dx \wedge d z}{\lambda f} - c\dfrac{z dx \wedge d y}{\lambda f} 
\]
up to scaling.
\end{defn}
\begin{table}[ht]

\begin{centering}
\begin{tabular}{c | c}
Polynomial equation $f$ & $(a,b,c)$ \\ [0.5ex] \hline  \\ [-2ex]
$x^3 + y^3 + z^3 + \tau xyz$  & $(1,1,1)$ \\[0.5ex]
$x^4 + y^4 + z^3 + \tau xyz $ & $(1,1,2)$ \\[0.5ex]
$x^6 + y^3 + z^3 + \tau xyz$ & $(1,2,3)$ \\[1ex]
\end{tabular}
\label{table:nonlin}
\end{centering}
\caption{Local models for elliptic singularities}
\end{table}
We will now characterize the Hodge theoretic properties of purely elliptic log symplectic pairs of dimension 4. Before proceeding, note that, by our definition, if $(X,Y)$ is a purely elliptic log symplectic pair, then the singular locus of $Y$ is a smooth curve.
\begin{theorem}\label{prop:pym}
Let $(X,Y)$ be a purely elliptic log symplectic pair with log symplectic form $\sigma$ and so that $\dim X = 4$. Let $C$ denote the curve of singularities in $Y$ and assume that $C$ is connected. Let $Y'$ be the proper transform of $Y$ under blow up $b : \Bl_CX \rightarrow X$ of $C$, and let $E$ denote the exceptional divisor of $b$. Then the following statements hold.
\begin{enumerate}
\item $(\Bl_CX, Y'\cup E)$ is an snc log symplectic pair. 
\item If $h^{2,0}(X) = 0$ then $(\Bl_CX,Y' \cup E)$ is log symplectic of pure weight 1.
\item $\Gr_3^W\HH^2(X \setminus Y;\mathbb{Q}) \cong \HH^1(C';\mathbb{Q})(-1)$ where $C'$ is an elliptic curve.
\end{enumerate}
\end{theorem}
\begin{proof}
For ease of notation, we will address only the case where $f= x^3 + y^3 + z^3 + \tau xyz$. The remaining cases are identical and we leave them to the reader.

Let us start by proving (1). We know that $\sigma$ has only nondegenerate logarithmic poles along $Y^\mathrm{sm}$ (the smooth part of $Y$) and is closed. Our goal is then to show that $b^*\sigma$ has only logarithmic poles along $Y'$ and $E$ and $b^*\sigma \wedge b^*\sigma$ is logarithmically nondegenerate.  We can check this locally, viewing $\mathrm{Bl}_C X$ as a subvariety of $\mathbb{C}^4 \times \mathbb{P}^3$ where $\mathbb{P}^3$ has coordinates $[q_1: q_2 : q_3: q_4]$, and $\mathrm{Bl}_C X$ has local equation
\[
xq_2 = yq_1, \quad xq_3 = zq_1, \quad yq_3 = zq_2, 
\]
and in fact only the first two equations are necessary. We can then cover this with charts where $q_1, q_2, q_3 = 1$ respectively. By symmetry, we will only look at $q_1 = 1$. Therefore, we have that $y = xq_2, z= xq_3$. Making this substitution, we find that 
\[
f(x, xq_2,xq_3) = x^3(1 + q_2^3 + q_3^3 + \tau q_2q_3),
\]
which is a smooth cubic. Let $f_0 = 1 + q_2^3 + q_3^3 + \tau q_2q_3$. Then $\sigma$ becomes
\[
(d \log f_0 + 3 d \log x) \wedge d w + d\log x \wedge ( (q_2 - q_3)(d q_2 - d q_3))/\lambda f_0 + (dq_2 \wedge d q_3) / \lambda f_0.
\]
Therefore $f_0xb^*\sigma$ is holomorphic. Furthermore, $d b^*\sigma$ is closed since $\sigma$ is closed, so, by definition, $b^*\sigma$ has logarithmic poles along $xf_0 = 0$.

Now we must show that $b^*\sigma \wedge b^*\sigma$ is nondegenerate. We note that 
\[
\sigma \wedge \sigma = -(x\partial_x f + y\partial_y f + z \partial_zf )(dx \wedge d y\wedge dz \wedge d w)/\lambda f^2.
\]
Since $f$ is a homogeneous polynomial of degree 3, $x\partial_x f + y \partial_y f + z \partial_z f = 3f$, therefore,
\[
\sigma \wedge \sigma = -3\dfrac{d x \wedge d y \wedge dz \wedge dw}{\lambda f}
\]
We then have that
\begin{align*}
b^*(\sigma \wedge \sigma)   = b^*\sigma \wedge b^*\sigma &= -9 \dfrac{d x \wedge ( x dq_2 + q_2 dx ) \wedge (x dq_3 + q_3 d x) \wedge d w}{\lambda x^3 f_0} \\
& = - 9\dfrac{dx \wedge dq_2 \wedge d q_3 \wedge d w}{\lambda x f_0}.
\end{align*}
Since $f_0 = 0$ is smooth, it follows that $b^*\sigma \wedge b^*\sigma$ is a logarithmic volume form near $x = f_0 = 0$ and along $x = 0$. This completes the proof since the holomorphic symplectic form $\sigma$ is nondegenerate away from singular points of $Y$.

Now we prove (2). We know that $Y_C = Y' \cup E$ and hence the maximal number of components in $Y_C$ which intersect one another is at most 2. Thus the spectral sequence computing the mixed Hodge structure on $\HH^*(\Bl_CX \setminus Y'\cup E;\mathbb{Q})$ (\cite[Theorem 8.34]{voisin}) shows that $\Gr^W_4\HH^2(\Bl_C X \setminus Y_C ; \mathbb{Q}) = 0$, hence $I^{2,2}(\HH^2(\Bl_CX \setminus Y_C;\mathbb{Q})) = 0$. If $h^{2,0}(X) = 0$, then the same is true for $\Bl_CX$ by \cite[Theorem 7.31]{voisin} and the fact that $C$ is a curve. Therefore, it follows that $I^{2,0}(\HH^2(\Bl_YX \setminus Y_C;\mathbb{C})) = 0$. Thus the only nonvanishing $I^{2,w}(\HH^2(\Bl_CX\setminus Y_C;\mathbb{Q}))$ is $I^{2,1}$ and therefore $(\Bl_CX,Y_C)$ is of pure weight 1.

Now let us prove (3). By the local description of the blow up along $C$ described above, it is clear that the intersection of $f_0 = 0$ and $x = 0$ (which we will denote $Z$) is isotrivially fibered over $C$ by an elliptic curve $C'$. Since $Y' \cup E$ is a normal crossings anticanonical divisor in $\Bl_CX$, it follows by adjunction that $Y' \cap E = Z$ is anticanonical in $Y'$ and $E$ respectively. Therefore, it is a Calabi--Yau surface, which means it is either a K3 surface or an abelian surface. Since no K3 surface admits a smooth elliptic fibration over a curve, it follows that $Z$ is an abelian surface and hence $C'$ is an elliptic curve. By construction, $E$ is a $\mathbb{P}^2$ bundle over $C$, hence $\HH^1(E;\mathbb{Q}) \cong \HH^3(E;\mathbb{Q})(1) \cong \HH^1(C;\mathbb{Q})$ and the Gysin homomorphism 
\begin{equation}\label{eq:tops}
\HH^1(Z;\mathbb{Q}) \cong \HH^1(C;\mathbb{Q}) \oplus \HH^1(C';\mathbb{Q}) \longrightarrow \HH^3(E;\mathbb{Q})
\end{equation}
is surjective. The standard spectral sequence computing the mixed Hodge structure on cohomology of $\Bl_CX \setminus Y_C$ (see e.g. \cite[Theorem 8.34]{voisin}) says that $\Gr^W_3\HH^2(\Bl_CX \setminus Y_C;\mathbb{Q})$ is isomorphic to the kernel of the direct sum of Gysin maps,
\[
\HH^1(Z;\mathbb{Q}) \cong \HH^1(C;\mathbb{Q}) \oplus \HH^1(C';\mathbb{Q}) \longrightarrow \HH^3(E;\mathbb{Q}) \oplus \HH^3(Y';\mathbb{Q}).
\]
Therefore, $\dim I^{2,1;2}, \dim I^{1,2;2}\leq 1$ by (\ref{eq:tops}). By Proposition \ref{prop:bounds} and (2), we know that 
\[
\dim I^{2,1;2}, \dim I^{1,2;2} \geq 1.
\]
Therefore, $\Gr^W_3\HH^2(\Bl_CX \setminus (Y' \cup E);\mathbb{Q})$ is isomorphic to $\HH^1(C';\mathbb{Q})$. This proves (3).

\end{proof}

 One of the goals of \cite{pym} is to restrict the structure of varieties which may admit purely elliptic log symplectic structures. The cohomological restrictions in Proposition \ref{prop:pw2} offer stringent restrictions on the structure of $(X,Y)$ and $C$. 

We also note that an example satisfying all of the assumptions in Theorem \ref{prop:pym} coming from work of Feigin and Odesski \cite{fo} is described in \cite{pym}.

\subsection{Constructing more examples}\label{sect:leafblow}

The goal of this section is to explain how one can produce new examples of log symplectic varieties from known examples. This implies that unlike the case of compact holomorphic symplectic manifolds, we can produce {\em many} families of log symplectic varieties once we are able to produce a single example. In particular, the examples in Section \ref{sect:purel} can be used to produce a host of new examples by appropriate blow up. 
It will be important to have a nice local form for the log symplectic form $\sigma$ near smooth points of $Y_i$. Such an expression is supplied by Goto.
\begin{lemma}[{Goto \cite[Lemma 1-2]{goto}}]\label{prop:darboux}
Let $(X,Y)$ be log symplectic with $\dim X = 2d$ and let $p$ be a smooth point in $Y$. Then there are holomorphic coordinates $x_1,\dots, x_{2d}$ near $p$ so that $Y= V(x_1)$ and
\[
\sigma = \dfrac{dx_1}{x_1} \wedge dx_2 + dx_3 \wedge dx_4 + \dots + dx_{2d-1} \wedge dx_{2d}.
\]
\end{lemma}
We will refer to coordinates in which $\sigma$ has the form of Lemma \ref{prop:darboux} as {\em log holomorphic Darboux coordinates}. If $(X,Y)$ is a log symplectic pair, then the log symplectic form $\sigma$ induces an isomorphism
\[
\mathscr{T}_X(-\log Y) \longrightarrow \Omega_X(\log Y)
\]
Therefore $\sigma$ is equivalent to a holomorphic bivector field $\sigma' \in \HH^0(X,\wedge^2\mathscr{T}_X(-\log Y))$, which produces the inverse of $\sigma$,
\[
\Omega_X(\log Y) \longrightarrow \mathscr{T}_X(-\log Y).
\]
We note that $\mathscr{T}_X(-\log Y)$ is a subsheaf of $\mathscr{T}_X$, therefore, this map can be extended to a homomorphism from $\Omega_X(\log Y)$ to $\mathscr{T}_X$ which is, of course, no longer an isomorphism.
\begin{defn}
Let $V$ be a compact subvariety of $X$, We say that $V$ is a closed symplectic leaf of $(X,Y)$ if at each point $p$ in a Zariski open dense subset of $V$, the image of the map
\[
\Omega_X(\log Y)_p \longrightarrow \mathscr{T}_{X,p}
\]
coincides with $i_*\mathscr{T}_{Z,p}$.
\end{defn}

Let $p \in Y$ be a smooth point and assume that $V$ is a symplectic leaf of $\sigma$. Then we can write $\sigma$ in holomorphic Darboux coordinates in a polydisc $\bm{D}$ near $p$ to see that the image of $\Omega_X(\log Y)_p \rightarrow \mathscr{T}_{X,p}$ is the span of $\partial x_2,\dots, \partial x_{2d}$. Therefore, $V$ is locally written as $V(z_1,z_2)$.

The first claim in the following theorem is a generalization of a result appearing in work of Hacking and Keel \cite{hk}.

\begin{theorem}\label{prop:hackingkeel}
Let $(X,Y)$ log symplectic pair with log symplectic form $\sigma$ and let $V$ be a closed symplectic leaf of $\sigma$ of codimension 2 in $X$, which is contained in $Y_1$ and which intersects all strata $Y_I \subseteq Y_1$ transversally. Let $b : \Bl_VX \rightarrow X$ be the blow up of $X$ in $V$, let $Y_V$ be the proper transform of $Y$ under $b$. Then:
\begin{enumerate}
\item $(\mathrm{Bl}_VX,Y_V)$ is log symplectic with symplectic form $b^*\sigma$.
\item If $(X,Y)$ is pure of weight $w$ then $(\Bl_VX,Y_{V})$ is also pure of weight $w$.
\end{enumerate}
\end{theorem}
\begin{proof}
First, we note that since $V$ intersects each $Y_I$ in $Y_1$ transversally, the proper transform of $Y$ under the blow up map remains snc. Now we show that $b^*\sigma$ is log symplectic.

We start by proving (1). This is essentially a local computation, so we restrict ourselves to local charts in $X$. Let us choose some open polydisc $D$ in $X$ centered at a smooth point $p$ in $Y$ which contains $Z$, and let us assume that we have chosen a Darboux coordinate system (Proposition \ref{prop:darboux}), so that
\[
\sigma = d\log x_1\wedge dx_2 + \dots  + dx_{2d-1} \wedge d x_{2d}.
\]
In this coordinate system we have that $Y_1 = V(x_1)$ and by the argument before the statement of the proposition, $V = V(x_1,x_2)$. Now we will blow up $X$ in $V$ and show that $b^*\sigma$ remains log symplectic. The blow up is covered by two charts
\[
U_1 = \{ x_1t = x_2 \} \in D \times \mathbb{C}_{t}, \qquad U_2 = \{ x_1 = x_2s \} \in D \times \mathbb{C}_{s}.
\]
In $U_1$, the proper transform of $x_1 = 0$ does not intersect the exceptional divisor, and in $U_2$, the proper transform of $x_1 = 0$ is given by $s = 0$. In $U_1$, $\pi^*_V\sigma$ is written as
\begin{align*}
b^*\sigma|_{U_1} &= d\log x_1 \wedge (x_1d t + t d x_1) + \dots + dx_{2d-1} \wedge d x_{2d} \\
& =  dx_1 \wedge d t + \dots + dx_{2d-1} \wedge d x_{2d}
\end{align*}
and in $U_2$, $b^*\sigma$ is written as
\[
b^*_V\sigma|_{U_2} = d\log s \wedge d x_2 + \dots + dx_{2d-1} \wedge d x_{2d}.
\]
Therefore, $b^*\sigma |_{C_1}$ and $b^*\sigma|_{C_2}$ are expressed in holomorphic and log holomophic Darboux coordinates respectively and thus $b^*\sigma$ is nondegenerate everywhere in the preimage of a Zariski open subset of $V$. Since near a point outside of $V$, $b$ is the identity, we also have that $b^*\sigma$ has log poles on $Y_V$ and is nondegenerate everywhere except perhaps at points which are in the preimage of the intersection of $V$ and components $Y_2,\dots, Y_k$, which is a subset of codimension 2.

Let $p \in \Bl_VX$. Then $b^*\sigma$ is nondegenerate at $p$ if and only if $b^*\sigma$ is a section of $\Omega^{2d}_X(\log Y)$ which does not vanish at $p$. By the argument above, the vanishing locus of $b^*\sigma$ is contained in the intersection of $V$ and $Y_2\cup \dots Y_k$, which has codimension at least 2 in $\Bl_V X$. Therefore, $b^*\sigma$ does not vanish at any point in $\Bl_VX$.

Now we prove (2). For the sake of simplicity, we will ignore coefficients in cohomology. Let $E_V$ be the intersection of $\Bl_VX \setminus Y_V$ and the exceptional divisor of $b$. We begin by noticing that $b|_{\Bl_VX \setminus (Y_V \cup E_V)}$ is a biholomorphic map with image $X \setminus Y$. We will use the notation $\iota$ to denote its inverse composed with the open embedding of $\Bl_VX \setminus (Y_V \cup E_V)$ into $\Bl_V X \setminus Y_V$.

We would like to show that the map $\iota^* : \HH^2(\Bl_VX \setminus Y_V) \rightarrow \HH^2(X \setminus Y)$ is injective on $F^2$, and that this map sends the class represented by $\iota^*\sigma$ to the class represented by $\sigma$. If we can show these two things, then the result comes from the following argument. The map $\iota^*$ induces a morphism of mixed Hodge structures. A morphism of mixed Hodge structures between $(V_1,F_1,W_2)$ and $(V_2, F_2,W_2)$ sends $I^{2,w}(V_1)$ to $I^{2,w}(V_2)$ and agrees with the induced morphism on Hodge and weight graded pieces, therefore, the induced maps $I^{2,i}(\HH^2(\Bl_VX \setminus Y)) \rightarrow I^{2,i}(\HH^2(X \setminus Y))$ are injective for each $i$. We can represent $[\iota^*\sigma]$ uniquely as a sum of classes $\sigma_{2,0} + \sigma_{2,1} + \sigma_{2,2}$ where $\sigma_{2,i} \in I^{2,i}(\HH^2(\Bl_V X \setminus Y_V))$. Since $\iota^*$ is injective on each summand and $[\sigma]$ is pure of weight $w$ for some $w \in \{0,1,2\}$, it follows that $\sigma_{2,i} = 0$ for $i \neq w$.

We start by showing that the pullback map $\HH^2(\Bl_VX \setminus Y_V) \rightarrow \HH^2(X\setminus Y)$ is injective on $F^2$, or equivalently, that the map $\HH^{2d-2}_c(X \setminus Y) \rightarrow \HH^{2d-2}_c(\Bl_VY\setminus Y_V)$ is surjective in $\Gr_F^{2d-4}$. By construction, we have that $X \setminus Y$ is a Zariski open subset of $\Bl_V Y \setminus Y_V$, and its complement is $E_V$, an $\mathbb{A}^1$ bundle over $V$. We have a long exact sequence of mixed Hodge structures (see \cite{fujiki} or \cite[pp. 138]{ps})
\[
\dots \longrightarrow \HH^i_c(X \setminus Y) \longrightarrow \HH^i_c(\Bl_VY \setminus Y_V) \longrightarrow \HH^i_c(E_V) \longrightarrow \dots 
\]
The mixed Hodge structure on $\HH^{2d-2}_c(E_V)$ is dual to that of $\HH^0(E_V)$, therefore, it is isomorphic to $\mathbb{Q}(1-d)$, hence $\Gr_F^{2d-4}\HH^2(E_V) = 0$, and it follows that the induced map $\Gr_F^{2d-4}\HH^{2d-2}_c(X \setminus Y) \rightarrow \Gr_F^{2d-4}\HH^{2d-2}(\Bl_VY \setminus Y_V)$ is surjective. By duality, we obtain the desired injectivity.  

Now $[\iota^*\sigma]$ induces a class in $F^2\HH^{2}(\Bl_VX \setminus Y_V)$, via the quasi-isomorphism
\[
\mathscr{A}^\bullet_{\Bl_VX}(\log Y_V) \hooklongrightarrow \mathscr{A}^\bullet_{\Bl_VX \setminus Y_V}
\]
obtained by pullback, and the fact that $\iota^*\sigma$ is closed. The map $\HH^i(\Bl_VX \setminus Y_V) \rightarrow \HH^i(X \setminus V)$ is a morphism of mixed Hodge structure and is induced, again, by pullback on forms,
\[
\mathscr{A}^\bullet_{\Bl_V X \setminus Y_V} \longrightarrow \mathscr{A}^\bullet_{X\setminus Y}
\]
so that class represented by $b^*\sigma$ in $\HH^2(\Bl_VX \setminus Y_V)$ is identified with the class of $\sigma$, since $\iota$ is an isomorphism away from the exceptional divisor $E_V$. This proves the result.

\end{proof}

\begin{example}[Blow ups of the Feigin--Odesski example]
Feigin and Odesski \cite{fo} constructed a collections log symplectic pairs and Pym \cite{pym} has shown that one of their examples is purely elliptic. Let $E$ be an elliptic curve and let $\phi: E \rightarrow \mathbb{P}^4$ be an embedding of $E$ determined by a line bundle of degree 5. Then we may let $Y = \mathrm{Sec}(E)$ be the secant variety of lines determined by $E$. Pym showed that there is a purely elliptic log symplectic structure on $\mathbb{P}^4$ whose degeneracy divisor is $\mathrm{Sec}(E)$.

Following Proposition \ref{prop:pym} we may blow up $\mathbb{P}^4$ in a curve to obtain a snc log symplectic pair whose boundary composed of a pair of divisors, both of which are fiber bundles over a product of elliptic curves. The symplectic leaves of this Poisson structure are $\mathbb{P}^1$-bundles over factors of this product of elliptic curves. We may then produce an infinite number of examples of log symplectic pairs of pure weight 1 and dimension 4 by the process described in this section. 
\end{example}
\begin{example}[Cluster varieties]\label{ex:cluster}
Let $N$ be lattice of rank $2d$, and let $\Sigma$ be a rational, complete fan inside of $M \otimes \mathbb{R}$, so that each cone $C$ of $\Sigma$ has generators which generate $\mathrm{span}_\mathbb{R}(C) \cap M$. Then the associated toric variety $X_\Sigma$ is smooth and projective. Let $\alpha$ be a holomorphic 2-form on $\mathbb{C}^{*2d}$ of the form
\[
\alpha = \sum_{1\leq i< j\leq 2d} a_{ij} d\log x_i\wedge d\log x_j
\]
and assume that the alternating matrix $A=[a_{ij}]$ is nondegenerate. It is known \cite[Example]{pym} that $\alpha$ extends to a log symplectic form on the pair $(X_\Sigma, Y_\Sigma)$ where $Y_\Sigma$ is the union of the torus orbits in $X_\Sigma$ of dimension $\leq 2d-1$.

Let $\bm{Z} = \{Z_1,\dots, Z_k\}$ be symplectic leaves of the Poisson structure on $X_\Sigma$ induced by $\alpha$, and assume that for each torus orbit, $T$ of $X_\Sigma$, and for all $i\neq j$, $Z_i \cap T$ and $Z_j\cap T$ intersect transversally. Then we then may repeatedly blow up along $Z_1,\dots, Z_j$ to obtain a new log symplectic pair $(X_{\Sigma, \bm{Z}},Y_{\Sigma,\bm{Z}})$. The variety $U_{\Sigma,\bm{Z}} :=X_{\Sigma,\bm{Z}} \setminus Y_{\Sigma,\bm{Z}}$ is called a {cluster variety} by Hacking and Keel \cite{hk}. The cohomology of such varieties will be studied in detail in \cite{ha2}.
\end{example}

\section{Degenerations of holomorphic symplectic manifolds}\label{sect:degenerations}

The goal of this section is to show that log symplectic pairs also occur as the limit mixed Hodge structure of good degenerations of holomorphic symplectic varieties. We will also show that the components of a good degeneration of holomorphic symplectic varieties are themselves log symplectic pairs of pure weight. This result should be of independent interest, as it characterizes the components of good degenerations of hyperk\"ahler varieties.

\subsection{Limit mixed Hodge structures}\label{Sect:lmhs}

We will now discuss breifly the construction of limit mixed Hodge structures and their cup products. More details will be provided in Appendix \ref{sect:mhc}.

\begin{defn}
A {\em degeneration} is a K\"ahler manifold $\mathscr{X}$ equipped with a proper complex analytic map $\pi: \mathscr{X} \rightarrow \Delta$ and so that $X_t = \pi^{-1}(t)$ is smooth and projective if $t \neq 0$. We say a degeneration is {\em semistable} if near each point in $X_0$, there is an open polydisc in $\mathscr{X}$ with variables $(z_1,\dots, z_d)$ in which $\pi = z_1 \dots z_k$ for some $k \leq d$.
\end{defn}
Associated to any semistable degeneration, there is a {\em limit mixed Hodge structure}. The precise definition of the limit mixed Hodge structure is given in Appendix \ref{sect:mhc}, but we will discuss the relevant details here. Let $\Omega^1_{\mathscr{X}/\Delta}(\log X_0) = \Omega^1_\mathscr{X}(\log X_0)/\pi^*\Omega_\Delta^1(\log 0)$, and define 
\[
\Omega_{\mathscr{X}/\Delta}^i(\log X_0) = \wedge^i \Omega^1_{\mathscr{X}/\Delta}(\log X_0), \quad \Lambda_{X_0}^i = \Omega_{\mathscr{X}/\Delta}^i(\log X_0) \otimes \mathscr{O}_{X_0}.
\]
The rational structure on the limit mixed Hodge structure is given by $\HH^j(X_\infty;\mathbb{Q})$ where $X_\infty$ is defined to be the coproduct $\mathrm{exp}^*(\mathscr{X}\setminus X_0)$ as in the following diagram
\[
\begin{tikzcd}
X_\infty \ar[r] \ar[d] & \mathscr{X}\setminus X_0 \ar[d,"\pi|_{\mathscr{X} \setminus X_0}"]\\
\mathfrak{h}\ar[r,"\mathrm{exp}"] & \Delta \setminus 0.
\end{tikzcd}
\]
The space $X_\infty$ is a locally trivial fibration over the complex upper half plane $\mathfrak{h}$, therefore its cohomology isomorphic to $\HH^j(X_t;\mathbb{Q})$ for a choice of fiber $X_t = \pi^{-1}(t)$ with $t \neq 0$. Then, by \cite[Proposition 2.16]{steenbrink}, 
\[
\mathbb{H}^j(X_0;\Lambda_{X_0}^\bullet) \cong \HH^j(X_\infty;\mathbb{Q}) \otimes \mathbb{C}
\]
for all $j$. Under this isomorphism, one may define the Hodge filtration on $\HH^j(X_\infty;\mathbb{C})$ by taking the stupid filtration on $\Lambda_{X_0}^\bullet$, i.e.
\[
F^p\Lambda_{X_0}^\bullet = (0\longrightarrow \dots\longrightarrow 0 \longrightarrow \Lambda^p_{X_0} \longrightarrow \Lambda^{p+1}_{X_0}\longrightarrow \dots)
\]
and letting 
\[
F^p\mathbb{H}^j(X_0;\Lambda_{X_0}^\bullet) = \mathrm{im}(\mathbb{H}^j(X_0; F^p\Lambda_{X_0}^\bullet) \longrightarrow \mathbb{H}^j(X_0; \Lambda_{X_0}^\bullet))
\]
The wedge product pairing on $\Omega_{\mathscr{X}}^\bullet(\log X_0)$ descends to a cup product pairing $\mu$ on $\Lambda_{X_0}^\bullet$, and it is straightforward that $F^{p_1}\Lambda_{X_0}^\bullet \otimes F^{p_2}\Lambda_{X_0}^\bullet$ maps to $F^{p_1+ p_2}\Lambda_{X_0}^\bullet$ under this pairing. Therefore, the cup product pairing induces maps
\begin{equation}\label{eq:limitmix}
\mathbb{H}^{j_1}(X_0;\Lambda_{X_0}^\bullet) \otimes \mathbb{H}^{j_2}(X_0;\Lambda_{X_0}^\bullet) \longrightarrow \mathbb{H}^{j_1+j_2}(X_0;\Lambda_{X_0}^\bullet)
\end{equation}
which preserves the Hodge filtration. We would like this to be a morphism of mixed Hodge structures, but showing this is not straightforward, since the weight filtration on the limit mixed Hodge structure is not naturally defined in terms of $\Lambda_{X_0}^\bullet$ (see \cite[pp.268]{ps}), but is defined terms of an auxiliary complex $\bm{s}(\mathscr{A}^{\bullet,\bullet})$ (see Section \ref{sect:mhc} for definition), and there does not seem to be a natural map  
\[
\bm{s}(\mathscr{A}^{\bullet,\bullet}) \otimes \bm{s}(\mathscr{A}^{\bullet,\bullet}) \longrightarrow \bm{s}(\mathscr{A}^{\bullet,\bullet})
\]
compatible with the cup product which preserves the weight filtration. This deficiency has been addressed by Fujisawa \cite{fuji,fuji2} (following El Zein \cite{elz}) by producing a weak cohomological mixed Hodge complex (by simplicial methods) for which there is a cup product map compatible with wedge product in $\Lambda_{X_0}^\bullet$. We refer to \cite[\S 2.4]{fuji2} for a precise description of this complex as well as the cup product map. We have the following result.
\begin{theorem}[Fujisawa, {\cite[Propositions 3.3, 4.4]{fuji2}}]\label{thm:fuji}
There is a weak cohomological mixed Hodge complex $\mathscr{K}$ whose underlying complex of $\mathbb{C}$-valued sheaves we denote $\mathscr{K}_\mathbb{C}$. There is a filtered quasi isomorphism of complexes,
\[
\Lambda_{X_0}^\bullet \longrightarrow \mathscr{K}_\mathbb{C}
\]
and a natural cup product maps
\[
\Phi : \mathscr{K}\otimes \mathscr{K} \longrightarrow \mathscr{K}
\]
which induce morphisms of mixed Hodge structures, and so that the induced map
\[
\HH^p(X_0, \Lambda_{X_0}^\bullet) \otimes \HH^q(X_0, \Lambda_{X_0}^\bullet) \longrightarrow \HH^{p+q}(X_0, \Lambda_{X_0}^\bullet)
\]
agrees with the map induced the wedge product $\mu$.
\end{theorem}
Consequently, for every $j_1,j_2$, there is a cup product map 
\[
\HH^{j_1}(X_\infty;\mathbb{Q}) \otimes \HH^{j_2}(X_\infty;\mathbb{Q}) \longrightarrow \HH^{j_1+j_2}(X_\infty;\mathbb{Q})
\]
which is a morphism of mixed Hodge structures, and so that, under the isomorphisms in (\ref{eq:limitmix}), this map is induced by $\mu$. We have the following corollary.
\begin{corollary}
Let $\pi: \mathscr{X}\rightarrow \Delta$ be a semistable degeneration. The graded group $\oplus_{i=0}^{2d}\HH^i(X_\infty;\mathbb{Q})$ admits the structure of a Hodge ring when each $\HH^i(X_\infty;\mathbb{Q})$ is equipped with the limit mixed Hodge structure and the cup product pairing $\mu$.
\end{corollary}
\begin{proof}
If $X_t$ is a smooth fiber of $\pi$ then we know that 
\[
\dim \Gr_F^p\HH^{p+q}(X_t;\mathbb{C}) = \dim \Gr_F^p\HH^{p+q}(X_\infty;\mathbb{C})
\]
for all $p,q$ by \cite[Theorem 2.18]{steenbrink}. Therefore, $\bm{H}_\mathscr{X}$ satisfies the first part of Definition \ref{defn:hdgrng}, since Definition \ref{defn:hdgrng}(1) holds for all smooth varieties. To check the second condition in Definition \ref{defn:hdgrng}, note that by Theorem \ref{thm:fuji}, cup product on $\bm{H}_\mathscr{X}$ induces morphisms of mixed Hodge structures.
\end{proof}

\subsection{Good degenerations and symplectic Hodge rings}

In this section, we will prove that if $\pi : \mathscr{X} \rightarrow \Delta$ is a good degeneration of holomorphic symplectic manifolds then the corresponding limit mixed Hodge structure is in fact a symplectic Hodge ring.

\begin{defn}[Nagai \cite{nag}]
Let $\mathscr{X} \rightarrow \Delta$ be a semistable degeneration whose smooth fibers are holomorphic symplectic manifolds of dimension $2d$. Let $\sigma \in \HH^0(X_0, \Omega_{\mathscr{X}/\Delta}^2(\log X_0))$. By abuse of notation, we will also use $\sigma$ to refer to the corresponding element of $\HH^0(X_0;\Lambda_{X_0}^\bullet)$. 
\begin{enumerate}
    \item A component $W$ of $X_0 = \pi^{-1}(0)$ is called a {\em good component} if at each point of $W$, $\sigma^d$ is nonvanishing.
    \item  We say that $\pi : \mathscr{X} \rightarrow \Delta$ is a {\em good degeneration} if all components of $X_0$ are good. In other words, $\sigma^d$ is nonvanishing on $X_0$.
\end{enumerate}
\end{defn}

\begin{remark}
We expect that good degenerations are rare in practice, but that it is easier to find good components of semistable degenerations. 
\end{remark}

Note that $\sigma$ is $d$-closed on $X_0$, therefore $\sigma$ determines a closed element of $\Gamma( \bm{s}(\mathcal{C}_{Gd}^\bullet\Lambda^\bullet_{X_0}))$. In fact, it determines a closed element of $\Gamma(\bm{s}(\mathcal{C}_{Gd}^\bullet\tau_{\geq 2}\Lambda^\bullet_{X_0}))$, and, in turn, in $F^2\mathbb{H}^2(X_0;\Lambda_{X_0})$. We say that $\sigma$ is of pure weight $w$ if its class $[\sigma]$ in $F^2\mathbb{H}^2(X_0;\Lambda_{X_0}^\bullet)$ is in $I^{2,w}$.

\begin{theorem}
Let $\pi :\mathscr{X} \rightarrow \Delta$ be a good degeneration. Then $\HH^*(X_\infty;\mathbb{Q})$ is a symplectic Hodge ring. If $\sigma$ is of pure weight $w$ then $\HH^*(X_\infty;\mathbb{Q})$ is a symplectic Hodge ring of pure weight $w$.
\end{theorem}
\begin{proof}
This proof that follows is nearly identical to Theorem \ref{thm:hodge-filt}. We must show that a symplectic form on $\mathscr{X}$ induces a symplectic element in this Hodge ring. The remainder of the proof is very similar to the proof of Theorem \ref{eq:filtis}, so we will only sketch it. Before proceeding, note that Theorem \ref{thm:fuji} also says that the wedge product pairing on complexes induces the complex part of a homomorphism of mixed Hodge complexes, $\HH^*(X_\infty;\mathbb{Q}) \otimes \HH^*(X_\infty;\mathbb{Q}) \rightarrow \HH^*(X_\infty;\mathbb{Q})$.
\begin{enumerate}
    \item First, we may see that $\wedge^k \sigma : \Lambda_{X_0}^{d-k} \rightarrow \Lambda_{X_0}^{d+k}$ is an isomorphism for all $k$. First, recall that $\Lambda_{X_0}^{d-k}$ is a locally free sheaf of rank ${2d \choose d-k}$. Nonvanishing of $\sigma^d$ means that the Pfaffian of $\sigma$ is nonvanishing at each point in $X_0$, therefore $\wedge^k \sigma$ is an isomorphism at each point, hence an isomorphism of sheaves.
    \item Next, we remark that the map of complexes $L_\sigma: \Lambda_{X_0}^\bullet \rightarrow \Lambda_{X_0}^{\bullet + 2}$ sending an element to its wedge product with $\sigma$ induces a homomorphism $\mathbb{H}^j(X_0;\Lambda_{X_0}^\bullet) \rightarrow \mathbb{H}^{j+2}(X_0;\Lambda_{X_0}^\bullet)$. This map is the cup product with the cohomology class $[\sigma]$. 
    \item Finally, one must show that the induced map 
    \[
    [\sigma]^{(i+j)} : \Gr_F^{d-i-j}\mathbb{H}^{d-i}(X_0;\Lambda_{X_0}^\bullet) \longrightarrow \Gr^{d+i+j}_F\mathbb{H}^{d+i+2j}(X_0;\Lambda_{X_0}^\bullet)
    \]
    is an isomorphism. This is identical to the proof of Theorem \ref{eq:filtis}..
\end{enumerate}
The fact that this Hodge ring is of pure weight if $\sigma$ is of pure weight is evident.
\end{proof}

\subsection{Degenerations of projective IHS manifolds}
In this section, we will specialize our results to degenerations of projective IHS manifolds. A {\em projective IHS manifold} is a smooth, simply connected, projective variety which is holomorphic symplectic and so that $\dim \HH^0(X;\Omega^i_X) = 0$ of $i$ is odd and $1$ if $i$ is even and $i \leq 2d$. Assume that $\pi: \mathscr{X} \rightarrow \Delta$ is a good degeneration and that $X_t$ is a projective IHS manifold, then $h^{2,0}(X_t)=1$ for a general $X_t$ of $\pi$. It follows from work of Steenbrink \cite[Theorem 2.18]{steenbrink} that $\dim \HH^0(X_t;\Omega_{X_t}^2) = \dim \HH^0(X_0;\Lambda_{X_0}^2)$, therefore $\dim \HH^0(X_0,\Lambda^2_{X_0}) \cong \mathbb{C}$, thus $\pi: \mathscr{X} \rightarrow \Delta$ is automatically of pure weight 0, 1, or $2$.

\begin{defn}
Let $\pi : \mathscr{X} \rightarrow \Delta$ be a semistable degeneration. Then the monodromy operator $T : \HH^i(X_\infty;\mathbb{Q}) \rightarrow \HH^i(X_\infty;\mathbb{Q})$ is unipotent; therefore $N_i = \log T_i$ is nilpotent. 
\begin{enumerate} 
\item The {\em order of nilpotency of $N_i$} is the smallest $j$ so that $N_i^j = 0$. We use $\nu_i$ denote the order of nilpotency of $N_i$. 
\item We say that $\pi:\mathscr{X}\rightarrow\Delta$ is of type I if $\nu_2 = 0$, type II if the $\nu_2 = 1$, and type III if $\nu_2 = 2$.
\end{enumerate}
\end{defn}
\begin{remark}
This terminology originates in work of Kulikov \cite{kul} where degenerations of K3 surfaces are classified.
\end{remark}
For any semistable degeneration $\pi: \mathscr{X} \rightarrow \Delta$ operator $N_m^k$ induces the weight filtration on $\HH^m(X_\infty;\mathbb{Q})$ (this is the main result of \cite{g-na}). Therefore, $N_m^k$ is an operator on $\HH^m(X_\infty;\mathbb{Q})$ which maps $W_i$ to $W_{i-2}$ and induces isomorphisms  
\begin{equation}\label{eq:charmon}
N_m^k : \Gr^W_{m+k}\HH^m(X_0;\mathbb{Q}) \rightarrow \Gr^W_{m-k}\HH^m(X_\infty;\mathbb{Q})
\end{equation}
for all $m,k$.
\begin{proposition}
If $\pi:\mathscr{X} \rightarrow \Delta$ is a good degeneration of projective IHS manifolds with symplectic form $\sigma$, then $\sigma$ is of pure weight $\nu_2$. 
\end{proposition}
\begin{proof}
Note that if $\nu_2 = 0$ then (\ref{eq:charmon}) implies that $\Gr^W_i\HH^2(X_\infty;\mathbb{Q}) = 0$ if $i \neq 2$. Therefore $\sigma$ must have pure weight 0. If $\nu_2 = 1$ then $\Gr^W_4\HH^2(X_\infty;\mathbb{Q}) = 0$ and $\Gr^W_3\HH^2(X_\infty;\mathbb{Q}) \neq 0$. Since $\dim \Gr^W_3\HH^2(X_\infty;\mathbb{Q}) = \dim I^{2,1;2} + \dim I^{1,1;2}$ (since $\Gr^W_3\HH^2(X_\infty;\mathbb{Q})$ carries a pure Hodge structure of weight 3, by Hodge symmetry, and by Proposition \ref{prop:Ipqprops}). Therefore, $I^{2,1;2} \cong F^2\HH^2(X_\infty;\mathbb{C})$, hence $\sigma$ is of pure weight 1. Finally, if $\nu_2 = 2$ then $\Gr^W_4\HH^2(X_\infty;\mathbb{Q}) \neq 0$ and $\dim \Gr^W_4\HH^2(X_\infty;\mathbb{Q}) = \dim I^{2,2;2}$ (since $\Gr^W_4\HH^2(X_\infty;\mathbb{Q})$ carries a pure Hodge structure of weight 4, by Hodge symmetry, and by Proposition \ref{prop:Ipqprops}). Therefore, $I^{2,2;2}\cong \mathbb{C}$ thus $\sigma$ is of pure weight 2.
\end{proof}

\begin{corollary}\label{cor:chl-deg}
If $\pi:\mathscr{X}\rightarrow \Delta$ is a type III degeneration of projective IHS manifolds, then the cohomology ring $\HH^*(X_\infty;\mathbb{Q})$ has the curious hard Lefschetz property.
\end{corollary}

\begin{remark}
This corollary can also be partially deduced from the main result of \cite{hlsy}.
\end{remark}
\subsection{Remarks on Nagai's conjecture}\label{sect:nagcon}

In \cite{nag}, Nagai made a conjecture about the structure of the monodromy weight filtration of degenerations of projective IHS manifolds. 
\begin{conjecture}[Nagai]\label{conj:nag}
Let $\pi:\mathscr{X} \rightarrow \Delta$ be a semistable degeneration of relative dimension $2d$ of projective IHS manifolds. Then $\nu_{2i} = i \nu_2$ for all $i \leq d$.
\end{conjecture}
In the case where $\nu_1 = 2$, this conjecture has two proofs that we know of. The first is work of Koll\'ar, Laza, Sacc\`a, and Voisin \cite{klsv}, and the second follows from a stronger result of Soldatenkov \cite{sold}. We may specialize Theorem \ref{thm:chl} to the case of good degenerations of projective IHS manifolds to deduce a third proof.
\begin{proposition}
Let $\pi : \mathscr{X} \rightarrow \Delta$ be a degeneration of holomorphic symplectic varieties of pure weight 2. Then $\nu_2 = 2$ and $\nu_{2i} = 2i$ for all $i \leq d$.
\end{proposition}
\begin{proof}
Since the weight filtration on $\HH^i(X_\infty;\mathbb{Q})$ can be identified with the monodromy weight filtration on a nearby fiber (\cite{g-na}), $\nu_{2i} + i$ is the same as the maximal $j$ for which $\Gr_{j}^W\HH^{2i}(X_\infty;\mathbb{Q})$ is nonzero. By Theorem \ref{thm:chl}, $\HH^{2i}(X_\infty;\mathbb{Q})$ has Hodge--Tate mixed Hodge structure, and $\sigma^i$ is in $I^{2i,2i;2i}$ for all $i \leq d$ by Proposition \ref{prop:mixedis}. Therefore, since $\Gr_F^{2i}\HH^{2i}(X_\infty;\mathbb{C}) \neq 0$ for all $i \leq d$ it follows that $\Gr^W_{4i}\HH^{2i}(X_\infty;\mathbb{Q}) \neq 0$ and $\nu_{2i} = 2i$.
\end{proof}

In the case where $\nu_1 = 1$, Nagai's conjecture has been proven for all known examples of projective IHS manifolds by Green, Kim, Laza, and Robles in \cite{gklr}. Theorems \ref{thm:ehl} and \ref{thm:chl}, when specialized to the case where the corresponding cohomology ring is the limit mixed Hodge structure of a good degeneration of projective IHS manifolds, provides some justification.
\begin{proposition}\label{thm:nag}
Let $\pi :\mathscr{X} \rightarrow \Delta$ be a good degeneration of projective IHS manifolds with $\nu_2 = 1$. Then for $0 \leq k \leq 2d$, $\nu_{2k} \leq \min\{d-1,k-1\}$.
\end{proposition}
\begin{proof}
Any projective IHS manifold of dimension $d$ has the property that $\Gr_F^{2k}\HH^{2k}(X;\mathbb{C}) \cong \mathbb{C}$ for all $0 \leq k \leq d$, and that $\dim \Gr_F^{2k+1} \HH^{2k+1}(X;\mathbb{C}) \cong 0$ for all $k$. We also know that a good type II degeneration of projective IHS manifolds produces a limit mixed Hodge structure which is a symplectic Hodge ring of pure weight 1. Therefore, we are in the situation of the second statement in Theorem \ref{thm:ehl} and it follows that if $d \leq k \leq 2d$ then $\Gr_j^W\HH^k(X_\infty;\mathbb{Q}) = 0$ if $j \geq k+d-1$. Since $\nu_k + k$ is equal to the highest possible $j$ so that $\Gr_j^W\HH^k(X_\infty;\mathbb{Q}) \neq 0$, the result follows.
\end{proof}
A simple consequence of this is then that Nagai's conjecture holds in certain degrees.
\begin{proposition}\label{prop:nagai}
Let $\pi : \mathscr{X} \rightarrow \Delta$ be a good degeneration of projective IHS manifolds with $\nu_{2} = 1$. Then $\nu_{2(d-1)} = (d-1)$, $\nu_{2d} = d$, and $\nu_4 = 2$. 
\end{proposition}
\begin{proof}
The statement that $\nu_{2d} = d$ follows from the fact that $\dim \Gr^W_{3d}\HH^{2d}(X_\infty;\mathbb{Q}) \neq 0$ by Proposition \ref{prop:bounds}, and that by Theorem \ref{thm:ehl}, $\dim \Gr^W_{j}\HH^{2d}(X_\infty;\mathbb{Q}) = 0$ if $j > 3d$. Similarly, by Proposition \ref{prop:bounds}, $\dim \Gr^W_{3(d-1)}\HH^{2(d-1)}(X_\infty;\mathbb{Q}) \neq 0$, and by Theorem \ref{thm:ehl}(2) shows that, under the conditions of the Corollary, $\Gr^W_j\HH^{2(d-1)}(X_\infty;\mathbb{Q}) = 0$ if $j> 3(d-1)$.

Finally, we show that $\nu_4 = 2$. Since $\dim \Gr^W_6\HH^4(X_\infty;\mathbb{Q}) \neq 0$, this is equivalent to showing that $\dim\Gr^W_7\HH^4(X_\infty;\mathbb{Q}) = \dim \Gr^W_8\HH^4(X_\infty;\mathbb{Q})= 0$. We know that 
\[
\dim \Gr^W_7\HH^4(X_\infty;\mathbb{Q}) = \sum_{p+q = 7} \dim I^{p,q;4}
\]
by Proposition \ref{prop:Ipqprops}. Since $\dim \Gr_F^p \HH^4(X_\infty;\mathbb{C}) = 0$ if $p > 4$, we must have $I^{p,q;4} \cong 0$ if $p > 4$ or if $q > 4$ by Hodge symmetry. By assumption, $\dim \Gr^4_F\HH^4(X_\infty;\mathbb{C}) = 1$ and by Proposition \ref{prop:Ipqprops},
\begin{equation}\label{eq:4flo}
\dim \Gr_F^4\HH^4(X_\infty;\mathbb{C}) = \sum_{r} I^{4,r;4}.
\end{equation}
By Proposition \ref{prop:bounds},  $\dim I^{4,2;4} = 1$. Therefore, by Proposition \ref{prop:Ipqprops}, we must have $\dim I^{4,3;4} = 0$. By Hodge symmetry, this means that $\dim I^{3,4;4} = 0$ as well.

Following the same line of reasoning as above, $\dim \Gr^W_8\HH^4(X_\infty;\mathbb{Q}) = \dim I^{4,4;4}$. Proposition \ref{prop:bounds} and Equation (\ref{eq:4flo}) then imply that $\dim I^{4,4;4} = 0$. Therefore, the proposition is proved.
\end{proof}
This is enough to conclude that Nagai's conjecture is true in low dimension.
\begin{corollary}
Let $\pi:\mathscr{X} \rightarrow \Delta$ be a good type II degeneration of projective IHS manifolds of relative dimension $\leq 8$. Then Nagai's conjecture is true.
\end{corollary}

\subsection{Pure weight and degenerations}
We  will now show that if a degeneration of holomorphic symplectic varieties is of pure weight $w$ then each good component of that degeneration is log symplectic of pure weight $w$. This result depends crucially on Theorem \ref{thm:isomhs}, which will be presented in Appendix \ref{sect:mhc}.

Let $\pi:\mathscr{X} \rightarrow \Delta$ be a snc degeneration, and let $\tau$ be a local section of $\Omega_{\mathscr{X}/\Delta}^2(\log X_0) = \Omega^2_{\mathscr{X}}(\log X_0)/\Omega_{\mathscr{X}}(\log X_0) \wedge \pi^*\Omega_\Delta(\log 0)$. Then 
\[
\wedge d\log\pi: \Omega^2_{\mathscr{X}/\Delta}(\log X_0)\longrightarrow \Omega^3_\mathscr{X}(\log X_0), \qquad \tau \longmapsto  \tau \wedge d \log \pi
\]
is a well-defined map from since $\pi^*\Omega_\Delta(\log 0)$ is spanned by $d\log \pi$. Let $W$ be a component of $X_0$ and let $\partial W$ be the intersection of $W$ with the singular locus of $X_0$. We can compose $\wedge d\log \pi$ with the residue map from $\Omega^3_{\mathscr{X}}(\log X_0)$ to $\Omega_W^2(\log \partial W)$ to obtain a morphism of sheaves
\[
\bm{r}_W : \Omega^2_{\mathscr{X}/\Delta}(\log X_0) \longrightarrow \Omega^2_{W}(\log \partial W).
\]

\begin{theorem}\label{thm:gooddegen}
Let $\pi : \mathscr{X} \rightarrow \Delta$ be a semistable degeneration with holomorphic symplectic form $\sigma$. Then for each good component $W$ of $X_0$, the pair $(W,\partial W)$ is snc log symplectic with log symplectic form $\bm{r}_W(\sigma)$. If $\sigma$ is of pure weight $w$ then so is $\bm{r}_W(\sigma)$.
\end{theorem}
\begin{proof}
Recall that in an analytic chart $U$ centered at a point $p \in X_0$ in which $X_0 = V(x_1\cdots x_k)$, the sheaf $\Lambda_{X_0}^2$ is spanned over $\mathscr{O}_U$ by 
\[
d\log x_1, \dots, d \log x_k, d x_{k+1} , \dots, d x_{2d}
\]
under the relation 
\[
d\log x_1 = - d \log x_2 - \dots - d\log x_k
\]
(see, e.g. \cite{steenbrink}). In other words, any form written locally is equivalent to one without a $d \log x_1$ factor. Let $W$ be given by $x_1 = 0$ in this chart. Then the form $\sigma \in \Lambda^2_{X_0}$ may be written locally at $p$ as the restriction of a form
\[
\widetilde{\sigma}  =  \sum_{2 \leq i< j \leq k}f_{i,j}(\bm{x})d \log x_i \wedge d\log x_j + \sum_{2 \leq i \leq k < j}f_{i,j}(\bm{x})d \log x_i \wedge d x_j + \sum_{k < i< j }f_{i,j}(\bm{x})d  x_i \wedge d x_j
\]
to $x_1\dots x_k = 0$. The fact that $\sigma$ is nondegenerate at $p$ as a section of $\Lambda_{X_0}^2$ means that the form
\[
\widetilde{\sigma}|_p  =  \sum_{2 \leq i< j \leq k}f_{i,j}(0)d \log x_i \wedge d\log x_j + \sum_{2 \leq i \leq k < j}f_{i,j}(0)d \log x_i \wedge d x_j + \sum_{k < i< j }f_{i,j}(0)d  x_i \wedge d x_j
\]
is nondegenerate as a form in 
\[
\bigwedge^2 \bigoplus_{i=2}^k \mathbb{C} (d\log x_i) \oplus \bigoplus_{j=k+1}^{2d} \mathbb{C} (d x_j)
\]
However, according to \cite[Lemma 3.1]{fr}, $\bm{r}_W(\sigma)$ is the restriction of $\widetilde{\sigma}$ to $x_0 = 0$. The holomorphic differentials in $\Omega_{W}^1(\log \partial W)$ are precisely $d\log x_2,\dots, d\log x_k, dx_{k+1}, \dots , dx_{2d}$, hence the condition that $\sigma$ is nondegenerate at $p$ implies that $\bm{r}_i(\sigma)$ is nondegenerate at $p$ as an element of $\Omega^2_{W}(\log \partial W)$.

Corollary \ref{cor:undermhs} then shows that if $[\omega]$ is of pure weight then so is $[\bm{r}_W(\sigma)]$. This completes the proof of the theorem.
\end{proof}

The following result intersects with results of \cite{klsv} in the case where a general fiber of $\pi : \mathscr{X}\rightarrow \Delta$ is a projective IHS manifold. When $X_t$ is not projective IHS this appears to be new.
\begin{corollary}[Koll\'ar--Laza--Sacc\`a--Voisin \cite{klsv}]
Let $\pi:\mathscr{X} \rightarrow \Delta$ be a good degeneration of holomorphic symplectic manifolds of dimension $2d$ and of pure weight $w$. Then the dimension of the dual intersection complex of the central fiber of $\pi$ is $dw$.
\end{corollary}
\begin{proof}
The dual intersection complex of $X_0$ has dimension one more than the dimension of the dual intersection complex $\partial W_i$. By Theorem \ref{thm:gooddegen}, each $(W_i,\partial W_i)$ are log symplectic of pure weight $m$. By Theorem \ref{thm:classification}, each $\partial W_i$ has dual intersection complex of dimension $dw-1$.
\end{proof}
\begin{remark}
The results of Koll\'ar, Laza, Sacc\`a, and Voisin \cite{klsv} are more general than ours in that they do not require that $X_0$ be simple normal crossings, only that $(\mathscr{X}, X_0)$ be a dlt pair. This condition is much more natural from the perspective of the minimal model program, and much more likely to be obtainable in practice. As it stands, only a few examples of good degenerations are known. In dimension 2 the famous Kulikov--Persson--Pinkham \cite{kul,perpink} classification of degenerations says that up to base change, all degenerations of K3 surfaces are birational to good degenerations. Furthermore, in \cite{nag}, Nagai has produced examples of good degenerations of Hilbert squares of K3 surfaces starting with type II degenerations of K3 surfaces.
\end{remark}

\appendix
\section{An auxiliary result}\label{sect:mhc}

In this section we prove a result which we were not able to find in the literature, but is a standard application of the machinery of cohomological mixed Hodge complexes. Cohomological mixed Hodge complexes are somewhat technical, and are not used in an essential way in the body of the paper, so in order to preserve the flow of exposition, this result is kept separate. Furthermore, it seems possible that this result is of independent interest, so this section has been written so that it can be read more or less independently of the main text.

\subsection{Statement of the main result}

Let $\pi:\mathscr{X}\rightarrow \Delta$ be a semistable degeneration. Then the limit mixed Hodge structure associated to $\pi$ (described in detail below) has complex part described by the hypercohomology of the complex $(\Lambda_{X_0}^\bullet,d)$. Recall that 
$\pi^*\Omega_\Delta^1(\log 0) \subset \Omega_\mathscr{X}^1(\log X_0)$, and that we have defined $\Omega_{\mathscr{X}/\Delta}^1(\log X_0):=\Omega_\mathscr{X}^1(\log X_0) / \pi^*\Omega^1_\Delta(\log 0)$. We then define
\begin{align*}
\Omega^p_{\mathscr{X}/\Delta}(\log X_0) & = \wedge^p\Omega^1_{\mathscr{X}/\Delta}(\log X_0) \cong \Omega_{\mathscr{X}/\Delta}^p(\log X_0)/(\pi^*\Omega_\Delta^1(\log 0) \wedge \Omega^{p-1}_{\mathscr{X}/\Delta}(\log X_0)) \\
\Lambda^p_{X_0} &= \Omega^p_{\mathscr{X}/\Delta}(\log X_0) \otimes \mathscr{O}_{X_0}. 
\end{align*}
For each irreducible component $W$ of $X_0$, we let $\partial W$ be the intersection of $W$ with the singular locus of $X_0$. We have a morphism of complexes,
\[
\bm{r}_W: \Lambda_{X_0}^\bullet \longrightarrow \Omega^\bullet_{W}(\log \partial W)
\]
which sends a form $\omega$ to $\mathrm{Res}_W(d\log\pi \wedge \omega)$. 
\begin{theorem}\label{thm:isomhs}
Let $\pi : \mathscr{X} \rightarrow \Delta$ be a semistable degeneration, let $W$ be an irreducible component of $X_0$. Then there is a morphism of mixed Hodge structures
\[
\bm{r}^\mathrm{MHS}_W: \HH^\ell(X_\infty;\mathbb{Q}) \longrightarrow \HH^\ell(W \setminus \partial W;\mathbb{Q})
\]
whose complexification is given by $\bm{r}_W$.
\end{theorem}

The proof of Theorem \ref{thm:isomhs} will be presented in Section \ref{sect:gd-lsp} after introducing the necessary machinery. 

\subsection{Cohomological mixed Hodge complexes}

We will now review some background on cohomological mixed Hodge complexes associated to normal crossing pairs. The reader may consult \cite[\S 3.3]{ps} or \cite[\S 3]{elz-ldt} for proof of the results in this section. 

The data of a cohomological ($\mathbb{Q}$-)mixed Hodge complex on a quasiprojective variety $X$ is; $\bm{K}_\mathbb{Q}$ in $\mathrm{D}^+(X;\mathbb{Q})$ equipped with an increasing weight filtration $W$; $\bm{K}_\mathbb{C}$ in $\mathrm{D}^+(X;\mathbb{C})$ equipped with an increasing weight filtration $W_\mathbb{C}$ and a decreasing Hodge filtration $F$; a quasiisomorphism of complexes $\alpha_{\bm{K}}: \bm{K}_\mathbb{Q}\otimes \mathbb{C} \rightarrow \bm{K}_\mathbb{C}$ which induces a quasiisomorphism between $W_i\bm{K}_\mathbb{Q}\otimes \mathbb{C}$ and $W_{\mathbb{C},i}\bm{K}_\mathbb{C}$. This data must also have the property that there is a pure Hodge structure on $\mathbb{H}^i(X,\Gr^W_j\bm{K}_\mathbb{Q})$ for all $i$ and $j$ induced by $F$. See \cite[Defintion 3.13]{ps} for a precise definition.

Given a cohomological mixed Hodge complex on $X$, the hypercohomology groups $\mathbb{H}^i(X,\bm{K}_\mathbb{Q})$ admit mixed Hodge structures.

Let $(\bm{K}_\mathbb{Q}, W, \bm{K}_\mathbb{C}, W_\mathbb{C}, F)$ and $(\bm{L}_\mathbb{Q}, W, \bm{L}_\mathbb{C}, W_\mathbb{C}, F)$ be a pair of mixed Hodge complexes. A morphism of mixed Hodge complexes is a pair of morphisms of complexes, $\phi_\mathbb{Q}: \bm{K}_\mathbb{Q} \rightarrow \bm{L}_\mathbb{Q}$ which preserves the filtration $W$, a morphism $\phi_\mathbb{C}: \bm{K}_\mathbb{C} \rightarrow \bm{L}_\mathbb{C}$ preserving $W_\mathbb{C}$ and $F$, so that the diagram
\[
\begin{tikzcd}
\bm{K}_\mathbb{Q} \otimes \mathbb{C} \ar[r,"\phi_\mathbb{Q} \otimes \mathbb{C}"] \ar[d,"\alpha_{\bm{K}}"] & \bm{L}_\mathbb{Q} \otimes \mathbb{C} \ar[d,"\alpha_{\bm{L}}"] \\ 
\bm{K}_\mathbb{C} \ar[r,"\phi_\mathbb{C}"] & \bm{L}_\mathbb{C}
\end{tikzcd}
\]
commutes up to homotopy (\cite[Definition 3.16]{ps}).

\subsection{Cohomology of a smooth quasiprojective variety}\label{sect:noncompmhs}
Let $U$ be a quasiprojective variety, and choose a simple normal crossings compactification $X$ of $U$. Let $Y = X \setminus U$ and let $j : U \hookrightarrow X$ be the injection map. Let $i: Y \hookrightarrow X$ be the closed embedding. The natural map
\[
\Omega^\bullet_X(\log Y)\hookrightarrow Rj_*\Omega^\bullet_U 
\]
is a quasiisomorphism, therefore $\mathbb{H}^i(X,\Omega_X^\bullet(\log Y)) \cong \HH^i(U;\mathbb{C})$ (\cite[Proposition 4.3]{ps}). We will use $\Omega^\bullet_X(\log Y)$ as the complex $\bm{K}_\mathbb{C}$ for the mixed Hodge complex that we wish to construct. Let $W_i^\mathbb{C}$ and $F^p$ be the filtrations on $\Omega_X^\bullet(\log Y)$ described in Definition \ref{defn:hodgeweight}

It is now necessary to produce a rational structure $\bm{K}_\mathbb{Q}$ with a weight filtration which matches $W^\mathbb{C}$. This can be done using a chain of quasi isomorphisms (as in \cite[\S 4.3]{ps}), but a more direct approach uses logarithmic structures (as in \cite[\S 4.4]{ps}). We will now explain this. We will let $\mathscr{M}_{X,Y}$ denote the sheaf of monoids given by $\mathscr{O}_X \cap j_*\mathscr{O}_U$, which in fact defines a log structure on $X$ called the divisorial log structure with respect to $Y$. We let $\mathscr{M}^\mathrm{gp}_{X,Y}$ denote the corresponding sheaf of abelian groups. We may identify $\mathscr{M}^\mathrm{gp}_{X,Y}$ with the sheaf of invertible local sections of $\mathscr{O}_X(*Y)$. There is a map
\[
\exp : \mathscr{O}_X \longrightarrow \mathscr{M}_{X,Y}^\mathrm{gp}, \quad f \longmapsto \exp(2\pi \mathtt{i} f)
\]
whose kernel is $\mathbb{Z}_X$ whose cokernel is $a_*\mathbb{Z}_{Y(1)}$. Here $Y(1)$ denotes the normalization of $Y$\footnote{Notation in this section differs from the previous sections slightly in order to be consistent with our main reference, \cite{ps}.} and $a: Y(1) \rightarrow X$ is the composition of the normalization map and the embedding of $Y$ into $X$. 
\begin{defn}
Let
\[
K_p^q = \mathrm{Sym}^{p-q}_\mathbb{Q}(\mathscr{O}_X) \otimes \bigwedge^q (\mathscr{M}_{X,Y}^\mathrm{gp} \otimes_\mathbb{Z}\mathbb{Q}).
\]
We obtain complexes $K_p^\bullet$ for each $p$ by taking the differential
\begin{equation}\label{eq:differentiallog}
d(f_1\dots f_{p-q} \otimes y) = \sum_{i=1}^{p-q} f_1\dots f_{i-1} f_{i+1} \dots f_{p-q} \otimes (\exp(f_i) \wedge y).
\end{equation}
\end{defn}
We have obvious inclusions of complexes,
\[
K_p^\bullet \hooklongrightarrow K_{p+1}^\bullet, \quad f_1\dots f_{p-q} \otimes y \mapsto 1\cdot f_1\dots f_{p-q} \otimes y.
\]
One sees that there are morphisms of complexes of sheaves given by
\[
\varphi_p : K_p^\bullet \longrightarrow \Omega_X^\bullet(\log Y), \quad \varphi_p(f_1\dots f_k \otimes y_1\wedge \dots \wedge y_p) = \dfrac{1}{(2\pi \mathtt{i})^q}\left(\prod_{i=1}^{p-q}f_i\right) \mathrm{d}\log y_1 \wedge \dots \mathrm{d} \log y_p.
\]
The image of $\varphi_p$ lies in $W^\mathbb{C}_{p}\Omega^\bullet_X(\log Y)$ and in fact (\cite[Theorem 4.15]{ps}) we obtain quasiisomorphisms
\[
K_p^\bullet \otimes \mathbb{C} \longrightarrow W^\mathbb{C}_{p}\Omega^\bullet_X(\log Y).
\]
We let $K_\infty^\bullet$ be the direct limit of the sheaves  $K_p^\bullet$ under these morphisms. If we let $W_pK_\infty^\bullet$ be the image of $K_p^\bullet$, we obtain a quasiisomorphisms of filtered complexes
\[
\phi_\infty : (K_\infty^\bullet, W) \otimes \mathbb{C} \longrightarrow (\Omega_X^\bullet(\log Y), W_\mathbb{C}).
\]
One can show (\cite[Lemma 4.6]{ps}) that $\Gr_m^W\Omega^\bullet_X(\log Y)$ is quasiisomorphic to $a_{m*}\Omega_{D(m)}^\bullet[-m]$, that this map preserves Hodge filtrations, and that the data
\[
\mathrm{Hdg}^\bullet_{X,Y}:=(K_\infty^\bullet, W, \Omega_X^\bullet(\log Y), W_\mathbb{C},F)
\]
is a mixed Hodge complex on $X$ whose underlying cohomology groups are $H^*(U;\mathbb{Q})$. 

\subsection{The cohomological mixed Hodge complex of a semistable degeneration}\label{sect:semistmhs}

Recall that in Section \ref{Sect:lmhs} we began to describe the cohomological mixed Hodge complex associated to a semistable degeneration. In this section, we will complete that description. The reader is asked to consult Section \ref{Sect:lmhs} for notation.  

We will first construct a complex which is quasiisomorphic to $\Lambda_{X_0}^\bullet$ which admits a Hodge and weight filtration. We begin by letting 
\[
W_{i}\Omega_{\mathscr{X}}^j(\log X_0) = \begin{cases} 0 & \mathrm{if } \quad m < 0 \\
\Omega^j_{\mathscr{X}}(\log X_0) & \mathrm{ if } \quad i \geq j \\
\Omega^{j-i}_{\mathscr{X}} \wedge \Omega^i_{\mathscr{X}}(\log X_0) & \mathrm{ if } \quad 0 \leq i \leq j\end{cases}
\]
We then define a double complex
\[
\mathscr{A}^{p,q} = \Omega_\mathscr{X}^{p+q+1}(\log X_0)/ W_p\Omega_\mathscr{X}^{p+q+1}(\log X_0), \quad p,q \geq 0
\]
with differentials
\[
d' : \mathscr{A}^{p,q} \longrightarrow \mathscr{A}^{p+1,q}, \quad d'' : \mathscr{A}^{p,q} \longrightarrow \mathscr{A}^{p,q+1}
\]
defined to be
\[
d'(\omega) = d\log \pi \wedge \omega, \qquad d''(\omega) = d\omega.
\]
We let $\bm{s}(\mathscr{A}^{\bullet,\bullet})$ denote the corresponding single complex of sheaves. We define a pair of filtrations on $\bm{s}(\mathscr{A}^{\bullet,\bullet})$. First, we let 
\[
W^\mathbb{C}_r\mathscr{A}^{p,q} = \mathrm{im}( W_{r+2p+1}\Omega_{\mathscr{X}}^{p+q+1}(\log X_0) \longrightarrow \mathscr{A}^{p,q})
\]
(this is called $W(M)$ in \cite{ps}), and we let 
\[
F^r\bm{s}(\mathscr{A}^{\bullet,\bullet}) = \oplus_{q \geq r} \mathscr{A}^{p,q}.
\]
We will let $\bm{L}_\mathbb{C}$ be $\bm{s}(\mathscr{A}^{\bullet,\bullet})$ equipped with the filtrations $W^\mathbb{C},F$. We first note that there are maps
\begin{equation}\label{eq:preres}
\mu : \Lambda_{X_0}^q \longrightarrow \mathscr{A}^{0,q}, \quad \omega \mapsto (-1)^q d\log \pi \wedge \omega
\end{equation}
which induces a quasiisomorphism of complexes from $\Lambda_{X_0}^\bullet$ to $\bm{s}(\mathscr{A}^{\bullet,\bullet})$ (\cite[pp. 269]{ps}). Moreover, equipping $\Lambda_{X_0}^\bullet$ with the truncation filtration
\[
F^r\Lambda_{X_0}^\bullet = (0 \longrightarrow \dots \longrightarrow 0 \longrightarrow \Lambda_{X_0}^r \longrightarrow \Lambda_{X_0}^{r+1} \longrightarrow  \dots),
\]
the map $\mu$ induces a filtered quasiisomorphism between $(\Lambda_{X_0}^\bullet, F)$ and $(\bm{s}(\mathscr{A}^{\bullet,\bullet}),F)$. We would now like to define $\bm{L}_\mathbb{Q}$ as in the case of an snc pair $(X,Y)$. We may proceed with the construction of $\mathscr{M}_{\mathscr{X},X_0}^\mathrm{gp}$ as before, and define
\[
L_p^q = \mathrm{Sym}_\mathbb{Q}^{p-q}(\mathscr{O}_{\mathscr{X}}) \otimes \bigwedge^q (\mathscr{M}_{\mathscr{X},X_0}^\mathrm{gp} \otimes_\mathbb{Z} \mathbb{Q}).
\]
Again, by equipping $L_p^q$ with a differential identical to that of (\ref{eq:differentiallog}), we obtain a complex $L_p^\bullet$ for all $p$, and by taking direct limits, we obtain a complex $L_\infty^\bullet$. We have that $\pi$ is, by definition, an element of $\mathscr{M}_{\mathscr{X},X_0}^\mathrm{gp}$, therefore, we obtain a global section of $L_1^1(1)$
\[
\widetilde{\theta} = 1 \otimes \pi \otimes 2\pi\mathtt{i}.
\]
(The notation $L_1^1(1)$ means $L_1^1$ tensored with $\mathbb{Q}(2\pi \mathtt{i})$). We then define complexes
\[
\mathscr{C}^{p,q} = (i^*L_\infty^{p+q+1}/ i^*L^{p+q+1}_p)(p+1), \, p\geq0, p+q \geq -1.
\]
We turn this into a double complex by adding the differentials,
\[
d' : \mathscr{C}^{p,q} \longrightarrow \mathscr{C}^{p+1,q}, \quad d'' : \mathscr{C}^{p,q} \longrightarrow \mathscr{C}^{p,q+1}
\]
defined as
\[
d'(x \otimes y) = x \otimes (\widetilde{\theta} \wedge y)
\]
and $d''$ is just the differential on $L_\infty^p$. Let $\bm{L}_\mathbb{Q} = \bm{s}(\mathscr{C}^{\bullet,\bullet})$ be the total complex with differential $\bm{d} = d' + d''$. We define a filtration on $\bm{L}_\mathbb{Q}$ by letting
\[
W_{r} \mathscr{C}^{p,q} = \mathrm{im}(i^*L^{p+q+1}_{r + 2p + 1} \longrightarrow \mathscr{C}^{p,q}).
\]
The compatibility between weight filtrations on $\Omega_{\mathscr{X}}^\bullet(\log X_0)$ and those on $L_\infty^p$, along with the quasiisomorphisms
\[
\phi_p: K_q^p \otimes \mathbb{C} \longrightarrow W^\mathbb{C}_q\Omega^p_\mathscr{X}(\log X_0)
\]
described in the previous section imply that there is a quasiisomorphism of filtered complexes,
\[
\phi: (\bm{L}_\mathbb{Q},W)\otimes \mathbb{C} \longrightarrow (\bm{L}_\mathbb{C}, W^\mathbb{C}).
\]
Furthermore, one can compute $\Gr^W_r\bm{L}_\mathbb{C}$ and $\Gr^{W^\mathbb{C}}_r\bm{L}_\mathbb{C}$ to see that their hypercohomology admits a pure Hodge structure with the induced Hodge filtration. Thus we have the following result.
\begin{theorem}[{\cite[Theorem 11.22]{ps}}]
Let $(\mathscr{X},\pi)$ be a semistable degeneration. Then
\[
\mathrm{Hdg}_{\mathscr{X},\pi}^\bullet = (\bm{L}_\mathbb{Q},W, \bm{L}_\mathbb{C}, W_\mathbb{C}, F)
\]
defines a cohomological mixed Hodge complex.
\end{theorem}

\subsection{A morphism of mixed Hodge complexes}\label{sect:gd-lsp}

In this section, we will prove Theorem \ref{thm:isomhs}. Let $W$ be an irreducible component of $X_0$ and we let $\partial W$ be the intersection of $W$ with the singular locus of $X_0$. Let $p$ be a point in $X_0$. We may write $\pi$ locally in a polydisc $\bm{D}$ centered at $p$ as $x_1\dots x_\ell$ for $ 0 \leq \ell \leq k$, therefore, in these coordinates, we can write $\Lambda_{X_0}^1$ as the $\mathscr{O}_{\bm{D}}$-linear span of 1-forms 
\[
 d \log x_1 \dots , d \log x_\ell,  d x_{\ell + 1} , \dots d x_{d}
\]
modulo the relation
\begin{equation}\label{ref:relation}
d \log x_1 + \dots + d \log x_\ell = 0.
\end{equation}
Now let us take a holomorphic 1-form $\widetilde{\tau} \in \Omega^1_{\mathscr{X}}(\log X_0)$ representing $\tau \in \Lambda^1_{X_0}$, given by
\[
\widetilde{\tau} = f_{1} d \log x_1  + \dots + f_\ell d \log x_\ell + f_{\ell+1} dx_{\ell+1} +\dots + d x_{d} 
\]
for some collection of holomorphic functions $f_{i}$. We may (locally) produce a holomorphic $2$-form on $Y_1 = V(x_1)$ by applying (\ref{ref:relation}) to replace $d\log x_1$ with $-d \log x_2 - \dots - d\log x_{2d}$, then setting $x_1= 0$ in each function $f_{i}$. Friedman \cite[Lemma 3.1]{fr} shows that this produces a well-defined map of sheaves,\footnote{Note that our notation differs from Friedman's slightly.}
\[
r^1_W : \Lambda^1_{X_0} \longrightarrow \Omega_{W}^1(\log \partial W).
\]
This may be extended to a morphism of complexes, $\bm{r}_W : \Lambda_{X_0}^\bullet \rightarrow \Omega_{W_1}^\bullet(\log \partial W_1)$. Friedman shows that the map $\bm{r}_W$ can be obtained as the composition of
\begin{equation}\label{wedgepi}
\wedge d\log \pi  : \Lambda^i_{X_0} \longrightarrow \Omega^{i+1}_{\mathscr{X}}(\log X_0), \quad \omega \longmapsto \omega \wedge d\log \pi 
\end{equation}
and
\begin{equation}\label{respi}
\mathrm{Res}_{W} : \Omega_\mathscr{X}^{i+1}(\log X_0) \longrightarrow \Omega_{W}^i(\log \partial W).
\end{equation}
We will now show that $\bm{r}_W$  underlies a morphisms of mixed Hodge structure. 

\begin{proof}[Proof of Theorem \ref{thm:isomhs}]
Precisely, we must show that, given the cohomological mixed Hodge complexes defined in Sections \ref{sect:semistmhs} and \ref{sect:noncompmhs},
\begin{enumerate}
    \item The map $\mathbf{Res}_W : \bm{s}(\mathscr{A}^{\bullet,\bullet}) \longrightarrow \Omega^{\bullet-1}_W(\log \partial W)$ preserves the Hodge filtration and weight filtration,
    \item There is a map $\mathbf{Res}^\mathbb{Q}_W : \bm{s}(\mathscr{C}^{\bullet,\bullet}) \longrightarrow K^{\bullet-1}_\infty$ which preserves the weight filtration,
    \item The diagram 
    \begin{equation}\label{mapmhs}
\begin{tikzcd}
\bm{s}(\mathscr{C}^{\bullet,\bullet})\otimes \mathbb{C} \ar[r,"\phi"] \ar[d,"\mathbf{Res}_{W}^\mathbb{Q}\otimes \mathbb{C}"] & \bm{s}(\mathscr{A}^{\bullet,\bullet}) \ar[d,"\mathbf{Res}_{W}"] \\K_\infty^{\bullet-1} \otimes \mathbb{C} \ar[r,"\varphi"] &\Omega_{W}^{\bullet-1}(\log \partial W)
\end{tikzcd}
\end{equation}
commutes.
\end{enumerate}
To begin, we show that the map 
\[
\mathbf{Res}_{W}: \bm{s}(\mathscr{A}^{\bullet,\bullet}) \longrightarrow \Omega^{\bullet - 1}_{W}(\log \partial W)
\]
sending
\[
\oplus_{p+q = j} \alpha_{p,q} \in \bigoplus_{p+q = j} \mathscr{A}^{p,q} \longmapsto \mathrm{Res}_{W}(\alpha_{0,j})
\]
is a morphism of $F$-filtered complexes. First, we check that this is a morphism of complexes. It is clear by definition that if $p \geq 1$ then $\mathbf{Res}_{W}(\bm{d}\alpha_{p,q}) = 0$. We also see that $\bm{d} \alpha_{0,q} = d\alpha_{0,q} \oplus (-1)^qd\log \pi \wedge \alpha_{0,q} \in \mathscr{A}^{0,q+1} \oplus \mathscr{A}^{1,q}$, so $\mathbf{Res}_{W_i}(\bm{d}\alpha_{0,q}) = \mathrm{Res}_{W_i}(d\alpha_{0,q})$. We now must show that
\begin{equation}\label{eq:diffcom}
d\mathrm{Res}_{W}(\alpha_{0,q}) = \mathrm{Res}_{W}(d\alpha_{0,q}).
\end{equation}
Choose a point $p \in X_0$ and local coordinates $(x_1,\dots, x_d)$ centered at a point $p \in W$ and so that in these coordinates, $W = V(x_1)$ and $X_0 = V(x_1\cdots x_k)$ for some $k \leq d$. Then if we have a $j$-form
\[
\tau = \sum_{\substack{ I \subset \{2,\dots, k\} \\ J \subset\{k+1,\dots, d\}\\ |I|+ |J| = j-1}} f_{I,J}(\bm{x}) d\log x_1 \wedge d\log x_I \wedge d x_J  + \sum_{\substack{ I \subset \{2,\dots, k\} \\ J \subset\{k+1,\dots, m\}\\ |I|+ |J| = j}} g_{I,J}(\bm{x})  d\log x_I \wedge d x_J.
\]
Then
\[
\mathrm{Res}_{W}(\tau) = \sum_{\substack{ I \subset \{2,\dots, k\} \\ J \subset\{k+1,\dots, d\}\\ |I|+ |J| = \ell-1}} f_{I,J}(0,x_2,\dots,x_d) d\log x_1 \wedge d\log x_I \wedge d x_J. 
\]
Here we have used the notation $d\log x_I$ to mean $\wedge_{i \in I} d\log x_i$ and similarly, $dx_J$ to mean $\wedge_{j \in J} dx_j$. Then (\ref{eq:diffcom}) is clear from the definition of the differential.Therefore $\mathbf{Res}_{W}$ is a morphism of complexes. 

It then follows directly from the definitions that $\mathbf{Res}_{W}$ respects the Hodge filtration; $F^r\bm{s}(\mathscr{A}^\bullet)^n = \oplus_{q \geq r, p +q =n}\mathscr{A}^{p,q}$, and note that $\mathscr{A}^{0,q} \cong \Omega^{q+1}_\mathscr{X}(\log X_0)$, so if $\alpha = \oplus_{q \geq r, p + q= n} \alpha_{p,q}$, then $\bm{r}_W(\alpha) = \mathrm{Res}_W(\alpha_{0,n})$ if $n \geq r$ and $0$ otherwise. So $\bm{r}_W(\alpha) \in F^r\Omega_W^n(\log \partial W)$.

Now let us define them morphism of $W$-filtered complexes alluded to above $\mathbf{Res}^\mathbb{Q}_{W} : \mathscr{C}^{p,q} \rightarrow K_\infty^{\bullet-1}$. To define the map $\mathbf{Res}^\mathbb{Q}_{W}$ we will imitate the definition of $\mathbf{Res}_{W}$. We assume that we are in a local chart on $\mathscr{X}$ in which the variables $y_1,\dots, y_d$ are such that $y_1$ is a local generator for the hypersurface $W$ and that in this chart, $y_1,\dots, y_k$ generate all components of $X_0$. An element of  $i^*K_\infty^{p+q+1}$ is represented by sums of local sections of the form
\begin{align*}
(f_1\dots f_\ell) \otimes y_{1}\wedge y_{i_2}\wedge \dots \wedge y_{i_{p+q+1}}, \quad &\{i_2,\dots, i_{p+q+1} \} \in \{2,\dots, k\}\\
(g_1\dots g_t) \otimes y_{j_1}\wedge \dots \wedge y_{j_{p+q+1}}, \quad &\{j_1,\dots, j_{p+q+1} \} \in \{2,\dots, k\},
\end{align*}
for $f_i,g_j$ germs of functions on $\mathscr{X}$ restricted to $X_0$. We may then define the residue map $\mathrm{Res}_{W_i}^\mathbb{Q}$ on an element of $i^*K_\infty^{p+q+1}$ by letting 
\begin{align*}
\mathbf{Res}_{W}((f_1\dots f_\ell) \otimes y_{1}\wedge y_{i_2}\wedge \dots \wedge y_{i_{p+q+1}}) &=\left( (f_1\dots f_\ell) \otimes y_{i_2}\wedge \dots \wedge y_{i_{p+q+1}}\right)|_{y_1 = 0} \\ 
\mathbf{Res}_{W}((g_1\dots g_t) \otimes y_{i_1}\wedge \dots \wedge y_{i_{p+q+1}}) & = 0.
\end{align*}
Clearly, $K_0^{p+q+1}$ is in the kernel of $\mathbf{Res}^\mathbb{Q}_{W_i}$, hence there is a well-defined map from $\mathscr{C}^{0,q}$ to $K_\infty^{p+q}$. Therefore, we may define a map
\[
\mathbf{Res}^\mathbb{Q}_{W} : \bm{s}(\mathscr{C}^{p,q}) \longrightarrow K_\infty^{\bullet}
\]
so that 
\[
\oplus_{p+q=\ell} c_{p,q} \in \bigoplus_{p+q = \ell} \mathscr{C}^{p,q} \mapsto \mathrm{Res}^\mathbb{Q}_{W}(c_{0,q}).
\]
The proof that this collection of maps extends to a morphism of complexes can be seen directly from the definitions and the local description of the map above.

Now we must check that this map preserves the weight filtration. Assume that $c= \oplus_{p+q = \ell} c_{p,q}$ is as in the previous equation and that each $c_{p,q} \in W_r\mathscr{C}^{p,q}$, which is to say that it can be lifted to a class
\[
i^*K^{\ell+1}_{r + 2p +1}
\]
up to Tate twist. By definition of $\mathbf{Res}_W^\mathbb{Q}$ we can ignore everything except $c_{0,\ell}$. In particular, $c_{0,\ell}$ is just an element of $i^*K^{\ell+1}_{r+1}$ which we recall means that it can be written as an element of
\[
\mathrm{Sym}_\mathbb{Q}^{\ell-r}(\mathscr{O}_\mathscr{X}) \otimes \bigwedge^{\ell+1}(\mathscr{M}^\mathrm{gp}_{\mathscr{X},X_0}\otimes_\mathbb{Z}\mathbb{Q}).
\]
So $\mathbf{Res}_W(c)$ is in the image of 
\[
\mathrm{Sym}^{\ell-r}_\mathbb{Q}(\mathscr{O}_\mathscr{X}) \bigwedge^{\ell}(\mathscr{M}^\mathrm{gp}_{W,\partial W} \otimes_\mathbb{Z}\mathbb{Q})
\]
in $K_\infty^\ell$, which is, by definition, $W_rK^\ell_\infty$. Therefore, the map $\mathbf{Res}^\mathbb{Q}_W$ respects the weight filtration. The proof that $\mathbf{Res}_W$ preserves the weight filtration on $\bm{s}(\mathscr{A}^{\bullet,\bullet})$ is nearly identical, so we will omit it. It is then straightforward to see (based on the definition of the maps $\phi_p$ and $\varphi_p$) that the diagram commutes.
\end{proof}
If $\omega$ is a closed global section of $\Lambda_{X_0}^n$ then it defines a closed element of $\bm{s}(\mathcal{C}_{Gd}(\bm{s}(\mathscr{A}^{\bullet,\bullet}))^n$ by the canonical quasiisomorphism
\[
\Lambda_{X_0}^\bullet \longrightarrow \bm{s}(\mathcal{C}_{Gd}(\Lambda_{X_0}^n)) \longrightarrow \bm{s}(\mathcal{C}_{Gd}(\bm{s}(\mathscr{A}^{\bullet,\bullet}))),
\]
and therefore an element of $\mathbb{H}^n(X_0;\mathbf{s}(\mathscr{A}^{\bullet,\bullet}))$ which we call $[\omega]$. The following is needed for the proof of Theorem \ref{thm:gooddegen}.
\begin{corollary}\label{cor:undermhs}
If $\omega$ is closed global section of $\Lambda_{X_0}^n$ for some $n$, and $[\omega] \in I^{n,j}(\HH^n(X_\infty;\mathbb{Q}))$, then $[\bm{r}_W(\omega)]$ is in $I^{n,j}(\HH^n(W\setminus \partial W;\mathbb{Q}))$.
\end{corollary}
\begin{proof}
Recall that the complex $\Lambda^\bullet_{X_0}$ is quasiisomorphic to $\bm{s}(\mathscr{A}^{\bullet,\bullet})$ via the map $\wedge d\log \pi$. Therefore, the composition of $\wedge d \log \pi$ and $\mathbf{Res}_{W_i}$ is precisely $\bm{r}_W$. 

Furthermore, the maps $\wedge d\log \pi: \Lambda_{X_0}^\bullet \rightarrow \bm{s}(\mathscr{A}^{\bullet,\bullet})$ and $\mathrm{Res}_W: \bm{s}(\mathscr{A}^{\bullet,\bullet}) \rightarrow \Omega^\bullet_W(\log \partial W)$ can be extended to maps $\mathcal{C}_{Gd}(\wedge d\log \pi)$  and $\mathcal{C}_{Gd}(\mathbf{Res}_W)$ between the corresponding Godement resolutions which make the diagram
\[
\begin{tikzcd}
\Lambda_{X_0}^\bullet \ar[r] \ar[d] & \bm{s}(\mathscr{A}^{\bullet,\bullet}) \ar[r] \ar[d] & \Omega_W^\bullet(\log \partial W) \ar[d] \\ 
\bm{s}(\mathcal{C}_{Gd}(\Lambda_{X_0}^\bullet)) \ar[r] & \bm{s}(\mathcal{C}_{Gd}(\bm{s}(\mathscr{A}^{\bullet,\bullet}))) \ar[r]  & \bm{s}(\mathcal{C}_{Gd}(\Omega_W^\bullet(\log \partial W))) 
\end{tikzcd}
\]
commute. Therefore After applying the functor $\HH^n(\Gamma(-))$ to the lower line of the equation above, we get the map 
\[
\bm{r}_W : \mathbb{H}^n(X_0;\Lambda^\bullet_{X_0}) \longrightarrow \mathbb{H}^n(X_0;\mathbf{s}(\mathscr{A}^{\bullet,\bullet})) \longrightarrow \mathbb{H}^n(W;\Omega_W^\bullet(\log \partial W)).
\]
Thus $[\bm{r}_W(\omega)] = \mathbf{Res}_W([\omega])$. By Theorem \ref{thm:isomhs}, the map $\mathbf{Res}_W$ is the complex part of a homomorphism of mixed Hodge structures, so if $[\omega] \in I^{n,j}(\HH^n(X_\infty;\mathbb{Q}))$ then Proposition \ref{prop:Ipqprops} shows that $[\bm{r}_W(\omega)]$ is an element of $I^{n,j}(\HH^n(W\setminus \partial W;\mathbb{Q})$.
\end{proof}

\bibliographystyle{alpha}
\bibliography{logsymp}

\end{document}